\documentclass[a4,12pt]{article}
\usepackage{amssymb}
\usepackage{latexsym}
\usepackage{epsfig}
\usepackage{graphicx}
\usepackage{psfrag}

\newcommand{\newtext}{}

\newcommand{\R}{{\mathbb R}}

\newcommand{\Sbb}{\mathbb S}
\newcommand{\cC}{{\cal C}}

\newcommand{\proofbox}{{$\Box$}}

\def\tilde{\widetilde}
\def \bfo {\begin {eqnarray*} }
\def \efo {\end {eqnarray*} }
\def \ba {\begin {eqnarray*} }
\def \ea {\end {eqnarray*} }
\def \beq {\begin {eqnarray}}
\def \eeq {\end {eqnarray}}
\def \supp {\hbox{supp}\,}

\def \dist {\hbox{dist}}

\def \det {\hbox{det}\,}

\def \e {\varepsilon}
\def \p {\partial}

\def\M{{\mathcal M}}

\def\stwo{\Sbb^2}
\def\sthree{\Sbb^3}

\newtheorem{definition}{Definition}[section]
\newtheorem{theorem}[definition]{Theorem}
\newtheorem{lemma}[definition]{Lemma}

\newtheorem{proposition}[definition]{Proposition}

\begin{document}
\title{Full-wave invisibility of active devices\\ at all frequencies}

\author{Allan Greenleaf\footnote{Department of Mathematics,
University of Rochester, Rochester, NY 14627, USA. Partially
supported by NSF grant
DMS-0551894. }\\
Yaroslav Kurylev\footnote{Department of Mathematical Sciences, Loughborough
University,
Loughborough, LE11 3TU, UK}
\\
Matti Lassas,\footnote{
Helsinki University of Technology, Institute of Mathematics,
P.O.Box 1100, FIN-02015, Finland. Partially supported by Academy of
Finland CoE Project 213476.}
\\
Gunther Uhlmann\footnote{Department of Mathematics,
University of Washington, Seattle, WA 98195, USA. Partially supported
by the NSF and a Walker
Family Endowed Professorship.}
}

\date{Final revision}

\maketitle

\begin{abstract}
There has recently been considerable interest in the possibility,
both theoretical and
practical, of invisibility (or ``cloaking") from observation by
electromagnetic (EM) waves.
Here, we prove invisibility with respect to solutions of the
Helmholtz and Maxwell's equations,
for several constructions of cloaking devices.
The basic idea, as in the papers \cite{GLU2, GLU3,Le,PSS1}, is to
use a singular transformation that pushes isotropic electromagnetic
parameters forward into singular, anisotropic ones. We
define the  notion of finite energy  solutions of the Helmholtz and Maxwell's
equations for such singular electromagnetic parameters, and
study the behavior of the solutions  on the entire domain, including the
cloaked region and its boundary.
{We show that, neglecting dispersion, the construction of
\cite{GLU3,PSS1} cloaks
passive objects, i.e., those without internal currents, at all
frequencies $k$.}
Due to the singularity of the metric, one needs to work with weak solutions.
Analyzing the behavior of such solutions inside the cloaked region, we
show that, depending on the chosen construction, there appear new
``hidden'' boundary conditions at the surface separating the cloaked
and uncloaked
regions. We also consider the effect  on invisibility of active
devices inside the cloaked  region, interpreted as
collections of sources and sinks or  internal currents.
When these conditions are overdetermined, as happens for Maxwell's
equations,  generic
internal currents prevent the existence of finite energy solutions and
invisibility is compromised.

We give two basic constructions for cloaking a
region $D$ contained in a
domain $\Omega\subset\mathbb R^n, n\ge 3$, from detection by measurements
made at $\p \Omega$  of Cauchy
data of waves on $\Omega$. These
constructions, the {\it single} and {\it double coatings}, correspond
to surrounding
either just the outer boundary
$\p D^+$ of the cloaked region, or both $\p D^+$ and $\p D^-$, with
metamaterials whose EM
material parameters (index of refraction or electric permittivity and
magnetic permeability)
are conformal to a singular Riemannian metric on $\Omega$.
   For the single coating
construction, invisibility holds for the Helmholtz equation, but
fails  for Maxwell's
equations with generic internal currents. However, invisibility can
be restored by
modifying the single coating construction, by either inserting a
physical surface at
$\p D^-$ or  using the double coating.  When cloaking an infinite cylinder,
invisibility results for Maxwell's equations are valid if the coating
material is
lined on $\p D^-$ with a surface satisfying the soft and  hard 
surface (SHS) boundary
condition, but in general not  without such a lining, even
for passive objects.

\end{abstract}

\section{Introduction}

There has recently been considerable interest \cite{AE,MN,Le,PSS1,MBW}
in the possibility, both theoretical and practical, of  a
region or object being shielded (or ``cloaked") from
detection via electromagnetic (EM) waves.
\cite[\S4]{GLU1} established a  non-tunnelling
result for time-independent Schr\"odinger operators with highly
singular potentials. This can be
interpreted  as cloaking, at any frequency, with  respect to
solutions of  the Helmholtz equation
using a layer of isotropic, negative index of refraction  material.
\cite{GLU2,GLU3} raised the possibility of passive objects being
undetectable, in the
context of electrical impedance tomography (EIT).  Motivated by
consideration of certain
degenerating families of Riemannian metrics, families of singular
conductivities, i.e., not
bounded below or above, were given and rigorous results  obtained
for the conductivity
equation of electrostatics, i.e., for waves of frequency zero.
A related example of a complete but noncompact two-dimensional 
Riemannian manifold
with boundary having the same Dirichlet-to-Neumann map as a compact 
one was given
in \cite{LTU}.

More recently, there has been exciting  work based on the availability of
metamaterials which allow fairly arbitrary behavior of EM material parameters.
The constructions in
\cite{Le} are based on conformal mapping in two dimensions and are
justified via change of
variables on the exterior of the cloaked region.
\cite{PSS1} also proposes a cloaking construction
for Maxwell's equations based on a singular transformation of the
original space,
again observing that, outside the cloaked region, the solutions of
the homogeneous
Maxwell equations in the original space become solutions of the transformed
equations. The transformation  used is the same as used previously in
\cite{GLU2, GLU3} in the
context of Calder\'on's inverse conductivity problem. The paper
\cite{PSS2} contains
analysis of cloaking
on the level of ray-tracing, while full wave numerical simulations
are discussed in \cite{CPSSP}.
Striking positive
experimental evidence for cloaking from microwaves has recently been
reported in \cite{SMJCPSS}.

Since
the metamaterials proposed to physically implement these constructions need
to be fabricated with a given wavelength, or narrow range of wavelengths, in
mind, it is natural to
consider this problem in the frequency domain.
   (As in the earlier works, dispersion,  i.e.,
dependance  of the EM material parameters on $k$, which is certainly present
for  metamaterials,  is neglected.)

The question we wish to consider is then whether,  at some
(or all) frequencies $k$, these constructions allow cloaking with
respect to solutions of the Helmholtz equation or
time-harmonic solutions of Maxwell's equations. {\newtext We  prove
that this indeed is the case
for the constructions of \cite{GLU3,PSS1}, as long as the object
being cloaked is passive; in
fact, for the Helmholtz equation, the object can be an active device
in the sense of having
sources and sinks. On the other hand, for Maxwell's equations with
generic internal currents,
invisibility in a physically realistic sense seems highly
problematic. We give several ways of
augmenting or modifying the original construction so as to obtain
invisibility for all internal
currents and at all frequencies.}

Due to the degeneracy of the equations at the surface of the cloaked region,
is important to
rigorously consider  weak solutions to the Helmholtz and Maxwell's
equations on
all of the domain, not just the exterior of the cloaked region. We
analyze  various
constructions for cloaking from observation on the level of
physically meaningful EM waves,
i.e., finite energy distributional solutions  of the  equations, showing  that
careful  formulation
of the problem  is   necessary both mathematically and for correct
understanding of the physical  phenomena.
It turns out that the cloaking structure imposes hidden boundary
conditions on such waves at the surface of the cloak.  When
these conditions are
overdetermined, finite energy solutions typically do not exist.
\newtext{ The time-domain physical interpretation of this is not 
entirely clear, but it seems to indicate an accumulation of energy or 
blow-up of the fields which would compromise the desired cloaking 
effect.}

As mentioned
earlier, the example in \cite{PSS1} turns out to be a special case of
one of the
constructions
from \cite{GLU2,GLU3}, which gave, in dimensions $n\ge 3$,
counterexamples to uniqueness for
Calder\'on's inverse problem \cite{C} for the conductivity equation.
(Such counterexamples have now also been
given for $n=2$ \cite{V,KSVW}.) Thus, since the equations of
electromagnetism (EM) reduce at frequency
$k=0$  to the conductivity equation with conductivity parameter
$\sigma(x)$, namely
$\nabla\cdot(\sigma\nabla u) =0$, for the electrical potential
$u$, the invisibility question  has already been
answered affirmatively in
the case of electrostatics.

The present work addresses the invisibility problem at all
frequencies $k\ne 0$.
We wish to cloak not just a passive object,
but rather an active ``device", interpreted as
a collection of sources and sinks, or an internal current, within $D$.
Since  the boundary value problems in
question may fail to have unique solutions (e.g., when $-k^2$ is a
Dirichlet eigenvalue on $D$), it is natural,  as in \cite{GLU1},
to use the set of
Cauchy data  at $\p\Omega$ of all of  the solutions, rather than the
Dirichlet-to-Neumann
operator on $\p\Omega$, which may not be well-defined.

The basic idea of \cite{GLU3,Le,PSS1} is to form new EM material parameters
by
pushing forward old ones
via a singular mapping $F$.
Solutions of the relevant  equations,
Helmholtz or Maxwell,
with respect to the old parameters then push forward to solutions of 
the modified
equations with respect to the new parameters outside the cloaked region.
However, when dealing with a {\it singular} mapping $F$, the converse is not
in general true.
This means that outside $D$, depending upon the
class of solutions considered,
there are solutions to the  equations
with respect to the new parameters
which are not the push forwards of  solutions to the  equations with
the old parameters.
   Furthermore, it is crucial that the solutions be dealt with on {\it
all} of $\Omega$,
and not only outside
   $D$.
Especially when dealing with the cloaking of active devices, this
gives rise to the question of  what are the proper {\it transmission}
conditions on $\p D$,  which  allow arbitrary  internal
sources to  be made invisible to an external observer.
To address these issues rigorously, one needs to make a suitable choice of the
class of {\it weak} solutions (on all of $\Omega$) to the
singular equations being considered.
For both mathematical and,  even more so,
physical reasons, the  weak solutions that are appropriate to
consider  seem to be  the {\it locally  finite energy} solutions; 
these belong to the Sobolev
space $H^1$ with respect to the singular volume form
$|\tilde g|^{1/2} dx$ on
$\Omega$ for Helmholtz, and $L^2(\Omega,|\tilde g|^{1/2} dx)$ for Maxwell.

These considerations are absent from \cite{Le,PSS1,PSS2}, where
the cloaking is justified by appealing to both the transformation of solutions
on the exterior of $D$ under smooth mappings $F$ (essentially the 
chain rule)  and the fact, in the high
frequency limit, that rays originating in $\Omega\setminus D$ avoid 
$\p D$ and do not enter $D$. As we
will show,  analysis of the  transmission conditions at
$\p D$ shows that the  constructions of \cite{PSS1,PSS2}, although adequate for
cloaking active devices in the absence of polarization, i.e., for 
Helmholtz, and
cloaking passive devices in the presence of polarization, i.e., for Maxwell,
fail to admit finite energy solutions
to Maxwell when generic active devices are present. Furthermore,
analysis of cloaking of an infinite cylinder, which was numerically explored in
\cite{CPSSP} and provides a model of the experimental verification of 
cloaking  in \cite{SMJCPSS},
shows that even cloaking a passive object may be problematic.
Fortunately, it is  possible to remedy the situation by augmenting or modifying
this construction.

We  now describe the results of this paper.  For what we call the
{\it single coating},
which is the construction of \cite{GLU3}, and apparently that
intended in \cite{PSS1}, we
establish invisibility with respect to the Helmholtz equation at all
frequencies $k\ne 0$.
In fact, one can not only cloak a
passive object in a region
$D\subset\subset\Omega$,  containing  material with   nonsingular index of
refraction $n(x)$, from all measurements  made at the
boundary $\p\Omega$,  but also an active device, interpreted as
a collection of sources and sinks within $D$.

Among the
phenomena described here is that finite energy solutions to the
single coating
construction must satisfy certain ``hidden'' boundary conditions at $\p D$. For
the Helmholtz equation, this is the Neumann boundary condition at $\p
D^-$, and it follows that
waves which propagate
inside $D$ and are incident to  $\p D^-$
behave as if the boundary were perfectly reflecting.
Thus, the cloaking structure on the exterior of $D$ produces a
``virtual surface" at $\p D^-$.
However, for Maxwell's equations with electric permittivity $\e(x)$
and magnetic
permeability $\mu(x)$, the situation is more complicated.
The hidden boundary condition forces the tangential components of both
the electric field {\bf{E}} and magnetic field {\bf{H}} to vanish on $\p D^-$.
For cloaking passive objects, for which the internal current $J=0$,
this condition can be satisfied, but for generic $J$, finite energy 
time-harmonic solutions
fail to exist, and
thus the single coating  construction is
insufficient for invisibility. In practice, even for cloaking passive objects,
this may degrade the effective invisibility.

We find two ways of dealing with this difficulty. 
One is to simply augment the single coating construction around a ball
   by adding a perfect electrical conductor (PEC)  lining at $\partial D$, in order to make the object inside the
   coating material to appear like a passive object.
(Such a lining was apparently incorporated, although claimed to be
unnecessary, into the code
used in \cite{CPSSP} in an effort to stabilize the numerics.)
However, for the sake of brevity, the necessary weak formulation of the
   boundary value problem for this setup will not
   be considered in this paper.
   
Alternatively,  one can
introduce a more elaborate construction, which we refer to as  the \emph{double
coating}.
Mathematically, this corresponds to a singular Riemannian metric which
degenerates in the same way as one approaches $\p D$ from both sides;
physically it would correspond to surrounding
\emph{both} the
inner and outer surfaces of $D$ with appropriately matched metamaterials.
We show that, for the double coating, no lining is necessary and  full
invisibility holds for arbitrary active devices, at all nonzero 
frequencies, for
both Helmholtz
and  Maxwell. It is even possible for the field to be identically
zero outside of $D$ while nonzero within $D$, and vice versa.

Finally, we also analyze cloaking within an infinitely
long cylinder, $D\subset\mathbb R^3$. In the main result of
\S\ref{cylinder} and
\S\ref{cylsec-shs},  we show that the cylinder $D$ becomes invisible
at all frequencies
if we use a double coating together with
the so-called {\it  soft and hard } (SHS) boundary condition
on $\p D$.
For the origin and properties of the SHS condition and a description of how
the SHS condition can be physically implemented, see \cite{HLS,Ki,Ki2,Li}.

We point out that there is some confusion in the physical
literature concerning the theoretical
possibility of invisibility. By this we mean uniqueness theorems for
the inverse problem of recovering the electromagnetic parameters from
boundary information (near field) or scattering (far field)
at a single frequency, or for all frequencies.
There is a vast literature
on this subject. We only mention here mathematical results directly
related  to the one mentioned
   in \cite{Le, SMJCPSS}.
The  Helmholtz operator at non-zero energy  for isotropic media is given by
$ \Delta+ k^2 n(x)$, where $n(x)$ is the index of refraction and
$k \ne 0$.  Unique
determination of $n(x)$ from boundary data  for a single frequency
$k$, under suitable
regularity assumptions on $n(x)$ and in
dimension $n\ge 3$, was proved in \cite{SU}, with a similar result
for the  acoustic wave
equation in \cite{N}.
(See
\cite{U} for a survey of  related results).
The article \cite{N} was referred to in \cite{Le,SMJCPSS}
as showing that perfect invisibility was not possible. However, the
results of \cite{SU,N} for the
Helmholtz equation are valid only under the assumption that the medium is
isotropic and that the index of refraction is bounded.
This does not contradict the possibility of invisibility
for an anisotropic index of refraction, nor for an unbounded isotropic
index of refraction. The constructions of \cite{GLU3,Le,PSS1} and
the present paper violate
both of these conditions. We also point out that the counterexamples 
given in \cite[~Sec. 4]{GLU1} yield invisibility for the
Helmholtz  equation, in dimension $n\ge 3$,
for certain  isotropic  {\it negative} indices of refraction which 
are highly singular (and negative)
on $\Omega\setminus D$.

We also note here that,
at fixed energy, the Cauchy data is equivalent to the inverse scattering data.
The connection between the fixed energy inverse scattering data, the
Dirichlet-to-Neumann map and the Cauchy data is discussed, for instance,  in
\cite{N} for the Schr\"odinger equation and in \cite{U} for the 
Helmholtz equation
in anisotropic media.  The scattering
operator is well defined for the degenerate metrics defined here;
see, e.g., \cite{M}.

There is a large literature (see \cite{U}) on uniqueness in the
Calder\'on problem for isotropic
conductivities under the assumption of positive upper and lower
bounds for $\sigma$. It was noted by
Luc Tartar (see
\cite{KV} for an account) that uniqueness fails badly if anisotropic
tensors are allowed, since if
$F:\overline\Omega\longrightarrow
\overline\Omega$ is a smooth diffeomorphism with $F|_{\p\Omega}= id$,
then $F_*\sigma$ and $\sigma$
have the same Dirichlet-to-Neumann map (and Cauchy data.) Note that since
$\e$ and $\mu$ transform in the same way,
this already constitutes a form of invisibility,
i.e., from the Cauchy data one cannot distinguish between
the EM material parameter pairs
$\e,  \mu$  and $\tilde \e= F_* \e,\, {\tilde \mu}=F_*\mu$.

Thus, uniqueness for anisotropic media, in the mathematical
literature, has come to mean uniqueness up
to a pushforward by a (sufficiently regular) map $F$.
Such uniqueness in the Calder\'on problem is known under various
regularity assumptions on the
anisotropic conductivity in
two dimensions \cite{S, N1, SuU,  ALP} and in three dimensions or higher
\cite{LU, LeU, LTU}, but for all of these results it is assumed that
the eigenvalues of $\sigma(x)$
are bounded below and above by positive constants.
Related to the Calder\'on problem is the Gel'fand problem, which uses Cauchy
data at all frequencies, rather than at a fixed one; for this 
problem, uniqueness
results are also available, e.g.,
\cite{BeK, KK}, with a detailed exposition in \cite{KKL}.  For example, in the
anisotropic inverse  conductivity problem as above,
Cauchy data at {\it all} frequencies determines the  tensor up to a 
diffeomorphism
$ F:\overline\Omega\longrightarrow
\overline\Omega$.

Thus, a
key point in the current works on invisibility that allows one to avoid the
known uniqueness theorems for the Calder\'on
problem is the lack of  positive lower and upper bounds on
the eigenvalues of these symmetric tensor fields.
In this paper, as in
\cite{GLU3,Le,PSS1}, the lower bound condition is violated
near $\p D$, and there fails to be a global diffeomorphism $F$
relating the pairs of material parameters having the
same Cauchy data.

For Maxwell's equations,
all of our constructions are made within the context of the
permittivity and permeability tensors $\e$ and
$\mu$ being conformal to each other, i.e.,
multiples of each other by a positive
scalar function;  this condition has been studied in detail in
\cite{KLS}. For Maxwell's equations in the time domain, this condition
corresponds to polarization-independent wave velocity. In particular,
all isotropic media are
included in this category.  This
seemingly special condition arises naturally from our construction,
since the pushforward $(\tilde\e,\tilde\mu)$ of an
isotropic pair $(\e,\mu)$ by a diffeomorphism need not be isotropic but
does satisfy this conformality.
For both mathematical and practical reasons, it would be very interesting to
understand cloaking for general anisotropic
materials in the absence of this assumption.

We believe that our results
suggest improvements which can be made in physical
implementations of cloaking.
In the very recent  experiment  \cite{SMJCPSS},
the configuration corresponds to a thin slice of  of an
infinite cylinder, inside of which a homogeneous, highly conducting disk
was placed in order to be cloaked.
This corresponds to the single coating with the metric $g_2$ (see \S2) on $D$
being a constant multiple of the Euclidian metric. The analysis here
suggests that lining
the inside surface
$\p D^-$ of the coating
with a material implementing the SHS boundary condition \cite{HLS,Ki,Ki2,Li}
should result in  less observable scattering
than  occurs without the SHS lining,
improving the partial invisibility  already observed.

The paper is organized as follows. In \S2 we describe the single and
double coating
constructions. We then establish  cloaking for the Helmholtz equation at all
frequencies in \S3. The notion of a finite energy solution
for the single coating is defined in \S\S3.2  and then the key step
for showing invisibility
is Proposition 3.5.  We discuss the Helmholtz equation for the double
coating In \S\S3.3;
there we define the notion of a weak solution and the  Neumann
boundary condition at the inner
surface of the cloaked region. The key step for invisibility for
Helmholtz at all frequencies
in the presence of the double coating is Proposition 3.11.

In \S4 we study invisibility at all frequencies for Maxwell's equations. We
define the notion of finite energy solutions for the single and
double coatings.
In \S5 we demonstrate invisibility for Maxwell's at all frequencies
for the double coating; see
Proposition 5.1. In \S6 we show that, for the single coating
construction, the Cauchy data
for Maxwell's equations must vanish on the surface of the cloaked region,
showing that generically finite energy solutions for Maxwell's
equations in the cloaked region
do not exist. In \S7 we consider an infinite cylindrical domain and
show invisibility
at all frequencies for Maxwell's equations for the double coating;
the key result is Proposition 7.1. In \S8, we consider how to cloak
the cylinder,
treating its surface as an obstacle with the SHS boundary condition.
Finally, in \S9, we briefly
indicate how general the constructions can be made. In particular, we
note that a modification
the double coating allows one to change the topology of the domain
and yet maintain
invisibility.

We would like to thank Bob Kohn for bringing the papers
\cite{Le,PSS1} to our
attention, and Ismo Lindell for discussions concerning the SHS
boundary  condition.

\section{Geometry and basic constructions}\label{gbc}

The material parameters of electromagnetism, namely the conductivity,
$\sigma(x)$;
electrical permittivity, $\e(x)$; and magnetic permeability, $\mu(x)$,
all transform as a product of a contravariant symmetric 2-tensor and
a $(+1)-$density.
That is, if $F:\Omega_1\longrightarrow\Omega_2,\quad y=F(x)$, is a
diffeomorphism between domains in $\R^n$, then $\sigma(x)=(\sigma^{jk}(x))$ on
$\Omega_1$ pushes forward to $(F_*\sigma)(y)$ on $\Omega_2$, given by

\beq\label{transf}
(F_*\sigma)^{jk}(y)=\left.
\frac 1{\det [\frac {\p F^j}{\p x^k}(x)]}
\sum_{p,q=1}^n \frac {\p F^j}{\p x^p}(x)
\,\frac {\p F^k}{\p x^q}(x)  \sigma^{pq}(x)\right|_{x=F^{-1}(y)},
\eeq
with the same transformation rule for the other material parameters.
It was observed by Luc
Tartar (se \cite{KV}) that it follows that if $F$ is a diffeomorphism
of a domain $\Omega$
fixing $\p\Omega$, then  $\sigma$ and $\tilde\sigma:= F_*\sigma$ have the same
Dirichlet-to-Neumann map, producing infinite-dimensional families of
indistinguishable
conductivities.

On the other hand, a Riemannian metric
$g=(g_{jk}(x))$ is a covariant symmetric two-tensor. Remarkably, in
dimension three or
higher, a material parameter tensor and a Riemannian metric can be associated
with each other  by
\beq\label{sigma-g}
\sigma^{jk}=|g|^{1/2}g^{jk},\quad\hbox{or}\quad
g^{jk}=|\sigma|^{2/(n-2)}\sigma^{jk},
\eeq
where $(g^{jk})=(g_{jk})^{-1}$ and $|g|=\det(g)$. Using this correspondence,
examples of singular anisotropic conductivities in $\R^n, n\ge 3$,
that are indistinguishable
from a constant isotropic conductivity, in that they have the same
Dirichlet-to-Neumann map, were given in \cite{GLU3}. The two constructions
there are based on  two different types
degenerations of Riemannian metrics, whose singular limits can be
considered as coming
from singular changes of variables. The singular conductivities
arising from these metrics via
the above correspondence are then indistinguishable from a constant
isotropic $\sigma$. In the
current paper, we will continue to examine one of these
constructions, correpsonnding to
pinching off a neck of a Riemannian manifold; we refer  to it as the
single coating.
We also introduce another construction, referred to as the double
coating. We start by
giving basic examples of each of these.
\bigskip

For both examples, let $\Omega=B(0,2)\subset \R^3$, the
ball of radius $2$ and center $0$, be the domain at the
boundary of which we make our observations;
$D=B(0,1)\subset\Omega$ the region to be cloaked; and
$\Sigma=\p D=\stwo$ the boundary of the cloaked region.
\bigskip

{\bf Single coating construction:} We begin by recalling an example from
\cite{GLU3,PSS1};  the two dimensional examples in
\cite{Le,V} are either essentially the same or closely related in structure.

For the single coating, we
blow up $0$ using the  map
\beq
\label{F1}
F_1:\overline{B}(0,2)\setminus\{0\}\to
\overline\Omega
\setminus \overline D,\quad
            \ F_1(x)=(\frac {r}2+1)\frac x{r},\quad r =|x|,\quad
0<r\le 2.
\eeq
On $\overline{B}(0,2)$, let $(g_e)_{ij}=\delta_{ij}$ be the Euclidian
metric, corresponding to constant isotropic material parameters; via the map
$F_1$,
$g_e$ pushes forward, i.e., pulls back by $F^{-1}$, to a metric on
$\overline{\Omega}\setminus\overline{D}$,
\ba
\tilde g_1=(F_1)_*g_e:=(F_1^{-1})^*(g_e)\quad .
\ea

Introducing the boundary normal coordinates $(\omega, \tau) $ in
the annulus $\Omega\setminus\overline{D}$, where $\omega=(\omega^1,
\omega^2)$ are local
coordinates on
$\Sigma= \stwo$ and $\tau >0$ is the distance in
metric $\tilde g_1$
to $\Sigma$,  we have, from (\ref{F1}),
\beq
\label{outside}
\tilde{g}_1=\tau^2\,d\omega^2+  d\tau^2,\quad
\tau =2(r-1).
\eeq
Here $d\omega^2=h_{\alpha \beta}(\omega) d\omega^\alpha
d\omega^\beta$ is the standard metric on
$\stwo$, induced by the  Euclidian metric on $\R^3$.
Note that $\tilde g_1$ has the following properties:

Consider a local $g_e$-orthonormal frame $(\p_r,v,w)$ on
$\overline{\Omega}\setminus\overline{D}$   consisting of the
radial vector
$$
\p_r=\frac \p{\p r}=\frac {x^j}{r}\frac \p{\p x^j}$$
and two vector fields
$v,w$. Then,
\beq
\label{propertyextra}
& &\tilde g_1(\p_r,\p_r)=4,\quad
\tilde g_1(\p_r,v)=\tilde g_1(\p_r,w)=0,\quad
\tilde g_1(w,v)=0,\\ \nonumber
& &\frac{\tilde g_1(v,v)}{(r-1)^2}\in [c_1,c_2],\quad
\frac{\tilde g_1(w,w)}{(r-1)^2}\in [c_1,c_2],
\eeq
where $c_1,c_2>0$. Thus, $\tilde g_1$ has one
eigenvalue bounded from below (with eigenvector corresponding
to the radial direction)
and two eigenvalues that are of order $(r-1)^2$
(with eigenspace $span \{v,w\}$).
In Euclidean coordinates,
we have that, for $|\tilde g_1|=\det(\tilde g_1)$,
\beq\label{property 1}
& &|\tilde g_1(x)|^{1/2}\sim C_1(r-1)^2,\\
& &|\tilde g_1^{ij}\nu_i|\leq C_2,\quad \nu_i=\frac{2x}{r}=2(\p_r)_i.\nonumber
\eeq
Here and below we use Einstein's summation convention, summing over
indices appearing both as sub- and super-indices in formulae,
and $\nu=(\nu_1,\nu_2,\nu_3)$ denotes the unit co-normal
vectors of  the surfaces $\{x\in \Omega\setminus \overline D :\ |x|=s\},$
$1<s<2$, with respect to the metric $\tilde g$.

On $D$, we simply let $\tilde g_2$ be the Euclidian metric. Together, the pair
$(\tilde g_1,\tilde g_2)$ define a singular Riemannian metric on
$\overline\Omega$,
\bfo
\tilde g=
\left\{
\begin{array}{cl}
     \tilde g_1, & x \in\Omega\setminus\overline{D}, \\
\tilde  g_2, & x \in D,
\end{array}
\right.
\efo
which is singular on
$\Sigma^+$, i.e., as
one approaches $\Sigma$ from $\Omega\setminus\overline{D}$;
in the sequel, we will identify the metric $\tilde g$
and the corresponding pair $(\tilde g_1,\tilde g_2)$.

     To unify  notation for the two basic constructions, we will denote
in the single coating case $M_1=\Omega$,
$M_2=D$ and let $M$ be the disjoint union $M=M_1\cup M_2$.
Also, for notational unity with the double coating, we let
$\gamma_1=\{0\}\subset M_1$,
$\gamma_2=\emptyset\subset M_2$, and $\gamma=\gamma_1\cup \gamma_2$.
Moreover, we denote
     $N_1=\Omega\setminus \overline D$,
$N_2=D$, $\Sigma=\p D$, and $N=N_1\cup \Sigma\cup N_2:=\Omega\subset \R^3$.

{\bf Double coating construction:}
The double coating refers to a  metric on $\Omega$ that
is degenerate on both sides of $\Sigma$ and has the same limit as one
approaches $\Sigma$ from both sides.

We now introduce  notation, shared with the single coating, that will
be used throughout
for the double coating. Let
$M_1=\Omega=B(0,2)$, which is compact with boundary, and
$M_2:=\sthree_{1/\pi}$, the 3-sphere of radius
$1/\pi$, which is compact without boundary, and
again let $M=M_1\cup M_2$ be their disjoint union. For the double coating,
let $\gamma_1=\{0\}\subset M_1$, $\gamma_2=\{NP\}\subset M_2$, where
$NP$ is a chosen point,
e.g.,  the North Pole of
$\sthree_{1/\pi}$, and $\gamma=\gamma_1\cup \gamma_2$.
As in the single coating example,
we let
$N_1=\Omega\setminus\overline{D}=B(0,2)\setminus\overline{B}(0,1)$,
$N_2=D=B(0,1)$, $\Sigma=\p D$, and $N=N_1\cup \Sigma\cup N_2\subset
\R^3$. We take the diffeomorphism $F_1:M_1\setminus
\gamma_1\longrightarrow
N_1$  to be the blow-up of $\gamma_1$ as
in the single coating, while we blow-up $\gamma_2$ by defining
$F_2:M_2\setminus
\gamma_2\longrightarrow N_2$ as follows.
Denote by $SP$ the point on $\Omega_2$ antipodal to $NP$.
Then the Riemannian normal coordinates centered at $SP$ are defined
on \linebreak$B(0,1)\subset T_{SP}\sthree\simeq \R^3$,
\bfo
\exp_{SP}:\,B(0,1) \to  M_2 \setminus\{NP\}.
\efo
Denote by $F_2$ the diffeomorphism
\bfo
F_2 = \left(\exp_{SP}\right)^{-1}: \, M_2 \setminus\{NP\} \to B(0,1).
\efo
Introduce (local)
spherical coordinates $ (\omega, r)$
on $N_2=B(0,1)$, considered as a subset of $T_{SP}(\sthree)$, with
$\omega =(\omega_1,
\omega_2), \,
\omega \in \Sigma= \p B(0,1)$ and $0 \leq r\leq 1$.
The standard metric $g$ on $\sthree_{1/\pi}$ in these coordinates
takes the form
\beq
\label{inside}
g_2=
\frac{\sin^2(\pi r)}{\pi ^2}\,
d\omega^2+dr^2,
\eeq
where $d\omega^2$ is again the standard metric on $\stwo$.

Observe that $\tilde g_2=(F_2)_*(g_2)$, as one approaches $\Sigma^-$
on $B(0,1)$,
has very similar properties to $\tilde g_1$ on $B(0,2)\setminus B(0,1)$ as one
approaches $\Sigma^+$. Indeed, again
consider the radial vector
$\p_r=\frac \p{\p r}=\frac {x^j}{r}\frac \p{\p x^j}$
at $x\in N_2$ and two vectors
$v,w$ such that in Euclidean metric $(\p_r,v,w)$
is a local orthonormal frame. Then,
as follows from (\ref{inside}),  at $x\in N_2$ with, say, $1/2 <r<1$,
\ba
& &\tilde g_2(\p_r,\p_r)=1,\quad
\tilde g_2(\p_r,v)=\tilde g_2(\p_r,w)=0,\\
& &\tilde g_2(w,v)=0,\quad
\frac{\tilde g_2(v,v)}{(1-r)^2},\,\,\frac{\tilde g_2(w,w)}{(1-r)^2}
\in [c_1,c_2],
\ea
where $c_1,c_2>0$. Thus, $\tilde g_2$ has one
eigenvalue bounded from below (with eigenvector corresponding
to the radial direction)
and two eigenvalues that are of order $(1-r)^2$,
and with respect to
the Euclidean coordinates on $N_2$,
\beq\label{property 1 b}
|\tilde g_2(x)|^{1/2}\leq C_1(1-r)^2,\quad
|\tilde g_2^{ij}\nu_i|\leq C_2,\quad \nu_i=-\frac{x_i}{r}=-(\p_r)_i,
\ \frac 12 <r<1.
\eeq

Set  $\tilde g_1=(F_1)_*g_e$ on $ N_1$, where $F_1$ is
defined as for the single coating example.
Together, these define a singular metric $\tilde
g=(\tilde{g}_1,\tilde{g}_2)$ on
the entire ball
$N=N_1\cup N_2 \cup \Sigma
=B(0,2)$.
Comparing (\ref{outside}) and (\ref{inside}), we see that, in the
Fermi coordinates \footnote{Recall that the Fermi coordinates associated
to $\Sigma$ are $ (\omega, \tau)$, where\ $\omega =(\omega^1, \omega^2)$
are local coordinates on $\Sigma $ and $ \tau=\tau(x)$ is the distance
from $x$ to $\Sigma$ with respect to the metric $\tilde g$,
multiplied by $+1$ in $N_1$ and $-1$ in $N_2$.} associated
to $\Sigma$, $|\tilde g|^{1/2} \tilde g^{ij}$ is
Lipschitz continuous on $N$;
note also that $|\tilde g|^{1/2}\tilde g_{ij}$ is not invertible at
$\p B(0,1)$.
\bigskip

Although they are distinct,  each of these constructions may be summarized as
follows. The domain $\Omega$, which we will refer to as $N$, decomposes as
$N=N_1\cup\Sigma\cup N_2$, where $N_1=\Omega\setminus \overline{D}, N_2=D$ and
$\Sigma=\p D$.
$N_1$ and $N_2$ are manifolds with boundary, with $\p
N_1=\p\Omega\cup \p D^+=\p N\cup \Sigma^+$ and $\p N_2=\Sigma^-$,
where the superscripts $\pm$ are used when  considering limits from
the exterior or interior of the cloaked region. The singular
electromagnetic material parameters on $N$ will correspond to a singular
Riemannian metric $\tilde g=(\tilde g_1,\tilde g_2)$, arising as the
pushforward
of a (nonsingular) Riemannian metric $g=(g_1,g_2)$
on a  manifold with two
components,
$M=M_1\cup M_2$, via a map $F:M\setminus\gamma\longrightarrow N$,
\bfo
F(x)=
\left\{
\begin{array}{cl}
     F_1(x), & x \in M_1\setminus \gamma, \\
     F_2(x), & x \in M_2\setminus \gamma.
\end{array}
\right.
\efo
Here, $M_1$ and
$M_2$ are disjoint, with $\overline{M_1}$ diffeomorphic to
$\overline{N}$; $\gamma_1=\gamma\cap M_1$ is either a point (the point being
blown up) for the single and double coatings, or a line (for the
cloaking of an infinite cylinder in \S\ref{cylinder},\ref{cylsec-shs}); and
$\gamma_2=\gamma\cap M_2$ is either empty (for the single coating) or a point
(for the double coating) or a line (for the cylinder.)
In \S\ref{appendix}, we will show that such
constructions exist in great generality, and for this reason the
proofs will be expressed in terms of analysis on $M$ and $N$.

In this generality, we say that $(M,N,F,\gamma,\Sigma,g)$ is
a \emph{coating construction} if $(M,g)$ is a (nonsingular)
Riemannian manifold, $\gamma\subset M$ and  $\Sigma\subset N$ are
as above, and
$F:M\setminus \gamma\to N\setminus \Sigma$ is diffeomorphism of either type.
This then
defines a singular Riemannian
metric $\tilde g$ everywhere on $N\setminus\Sigma=N_1 \cup N_2$, by
\bfo
\tilde g=
\left\{
\begin{array}{cl}
\tilde g_1:={F_1}_*g_1, & x \in N_1, \\
\tilde g_2:={F_2}_* g_2, & x \in N_2.
\end{array}
\right.
\efo
{If we introduce Fermi coordinates $(\omega,\tau)$ near $\Sigma$ as
above, the $\tilde g$ satisfies (\ref{propertyextra}),(\ref{property
1}) or (\ref{property 1 b}), with $r-1$ replaced by $\tau$, for the
single and double coatings, resp.
    From these, one sees that
$|\tilde g|^{1/2} g^{jk}$  has a jump discontinuity  across
$\Sigma$ for the single coating and is Lipschitz for the double coating.} Note
that in both examples, $N=
\Omega= B({0},2)$, so that $N$ and $M_1$ have the same topology.
However, in a direct extension of the double coating construction,
described in  \S\ref{appendix}, the domain $N$ containing the cloaked
region $N_2$ need not even
be diffeomorphic to $M_1$.

\section{The Helmholtz equation}

We are interested in invisibility of a cloaked
region with respect to the Cauchy data of solutions of the Helmholtz equation,
\beq\label{helm-g}
& &(\Delta_{g}+k^2) u=f\quad\hbox{in }\Omega,
\eeq
where $f$ represents a collection of sources and sinks. The
\emph{Cauchy data} $\cC_{g,f}^k$
consists of the set of pairs of boundary
measurements $\left(u|_{\p\Omega},\p_\nu  u|_{\p\Omega}\right)$ where $u$
ranges over solutions to (\ref{helm-g}) in some function or distribution space
(discussed below).
{
Let $(M,N,F,\gamma,\Sigma,g)$ be a single coating construction as in
\S\ref{gbc}.
For the moment, as
in the Introduction, we continue to refer to $N$ as $\Omega$, $N_2$
as $D$ and $\Sigma^+$
as $\p D^+$; we may assume that $M_1=N$, $M_2=D$ and $F_2= id$, so that
$\tilde{g}_2=g_2$ is a
(nonsingular) Riemannian metric on $D$. Thus,
$\tilde g$ is a
Riemannian metric on $\Omega$, singular on $\Omega\setminus D$, resulting from
blowing up the metric $g_1$ on $\Omega$ with respect to a point ${O}$
and inserting the $(D,g_2)$ into the resulting  ``hole".

We wish to show that $\cC_{\tilde g,\tilde f}^k=\cC_{g,0}^k$
for all frequencies $0<k<\infty$, if
$\hbox{supp}( \tilde f) \subset D$ and $k$ is not
a Neumann eigenvalue of  $(D,g_2)$.
Due to the singularity of $\tilde g$, it is necessary to consider
nonclassical solutions to
(\ref{helm-g}), and we will see that the exact notion of weak
solution is crucial. Furthermore, a hidden Neumann boundary condition on
$\p D^-$ is required for the existence of finite energy solutions.
Physically,
this means that  the coating on
$\Omega\setminus\overline{D}$ makes  the inner boundary
$\p D^-$ appear to be a perfectly reflecting ``sound-hard surface'' for waves
propagating in $D$, while, from the  exterior, the cloaked device is
invisible; that is, measurements of solutions of the Helmholtz equation at
$\p\Omega$ cannot distinguish between
$(\Omega,\tilde g)$ and $(\Omega, g)$.

\subsection{$k=0$ and weak solutions}

First consider the case when $k=0$ and $f=0$. As described in the
Introduction, this situation was treated in
\cite{GLU3} in the context of electrical impedance tomography.
There, it sufficed
to consider as weak solutions
those $L^\infty$ functions satisfying (\ref{helm-g}) (for the metric $\tilde
g$) in the sense of
distributions. It was shown that, for given Dirichlet data $h$ on
$\p\Omega$, (\ref{helm-g})
has a unique such solution, $\tilde u$, which must, by removable
singularity considerations, be constant on
$D$. These same conclusions would have held if we had considered the
larger class of spatial
$H^1$ weak solutions (defined below). However, for $k>0$ or $f\ne 0$, we
will see that this notion of
weak solution is inappropriate.

\subsubsection{  $k>0$ and spatial $H^1$ solutions}

\begin{definition}\label{spatial-def}
$\tilde u$ is a  \emph{spatial} $H^1$ solution to the Dirichlet
problem for the Helmholtz
equation,
\beq\label{helm-tg}
(\Delta_{\tilde g}+k^2)\tilde u=\tilde f\quad\hbox{on
}\Omega,\,
\tilde u|_{\p \Omega}=h
\eeq
if
\beq\label{h1-spa}
\tilde u\in H^1(\Omega,dx)
\eeq
and
\beq\label{distA}
\p_i(|\tilde g|^{1/2}\tilde g^{ij}\p_j\tilde u)
+k^2|\tilde g|^{1/2}\tilde u=|\tilde g|^{1/2} \tilde f\quad\hbox{in
}H^{-1}(\Omega, dx).
\eeq
\end{definition}

Here, for $s\in\mathbb R$, $H^s(\Omega,dx)$ refers to the
standard  Sobolev space of
distributions with $s$ derivatives in $L^2(\Omega, dx)$.
Note that (\ref{h1-spa}), together  with the properties of the metric tensor
given in \S\ref{gbc},
implies that
$|\tilde g|^{1/2}\tilde
g^{ij}\p_i\tilde u \in
L^2(\Omega, dx)$.

Later in our analysis (see  (\ref{Neumann boundary condition
on M_2})), we will see that (\ref{distA})
implies
that the normal derivative of $\tilde u$  on $\p D^-$
vanishes,
\ba
\p_r \tilde u|_{\p D^-}=0.
\ea
On the other hand,
the fact that $\tilde u\in H^1(\Omega,dx)$ implies that
\ba
\tilde u|_{\p D^-}= \tilde u|_{\p D^+}=\hbox{constant}:=u(O),
\ea
with $u$  the solution to $(\Delta_g+k^2)u=0$ in $\p \Omega$,
$u|_{\p \Omega}=h$, where the first equality follows from
the trace theorem for $H^1$ functions and the second from
considerations similar to those in \cite[Prop. 1]{GLU3}.
Note that, for generic $k$ and $h$, $u(O) \neq 0$.
Thus, $\tilde u_2:=\tilde u|_{D}$ needs to be a solution of the overdetermined
elliptic boundary value problem on $(D,\tilde{g}_2)$,
\beq
\label{overdetermined}
(\Delta +k^2)\tilde u_2=0, \quad \p_\nu \tilde u_2|_{\p D}=0,\,\,
\tilde u_2|_{\p D} =\hbox{constant} \neq 0.
\eeq
Clearly,  for  generic $k>0$ there exists no solution to
(\ref{overdetermined}) and therefore there
is no   weak solution to (\ref{helm-tg})  in
the sense of Definition
\ref{spatial-def}. Rather, one needs to use an $H^1$ norm adapted to
the singular Riemannian metric
$\tilde g$; this is in fact physically  natural, being essentially the
energy of the wave. We  formulate the correct notion in the next section.

\subsection{Finite energy solutions for the single coating}

We now give a more satisfactory definition of weak
solution, restricting the
notion to those solutions that are physically meaningful in that they
have  finite energy.

We now revert to the notation of $M,N,\dots$ when discussing the single coating
construction, i.e., let $(M,N,F,\gamma,\Sigma,g)$ denote a single coating as
in \S\ref{gbc}. Our first task is to understand in what sense the  expression
$|\tilde g|^{1/2}\tilde g^{ij}\p_i\tilde u$ is rigorously
defined.

To this end, define for ${\tilde \phi}\in  C^\infty(\overline N)$
\ba
\|{\tilde \phi} \|_X^2:=\int_{N}
(|\tilde g|^{1/2}\tilde g^{ij}\p_i{\tilde \phi} \p_j{\tilde \phi}+
|\tilde g|^{1/2}|{\tilde \phi}|^2)\,dx.
\ea
Let
\ba
H^1(N,|\tilde g|^{1/2}dx)=X:=\hbox{cl}_{X}( C^\infty(\overline N))
\ea
be the completion of  $ C^\infty(\overline N)$ with respect to the norm
$\|\,\cdotp\|_X$. We note that $H^1(N,|\tilde g|^{1/2}dx)\subset
L^2(N,|\tilde g|^{1/2}dx)$, so we can consider its elements
as measurable functions on $N$.

\begin{lemma}\label{lem: extension}
The map
\ba
\phi\longrightarrow D_{\tilde g}{\tilde \phi}=(D^j_{\tilde g}{\tilde
\phi})_{j=1}^3=
(|\tilde g|^{1/2}\tilde g^{ij}\p_i{\tilde \phi})_{j=1}^3,\quad \phi\in
C^\infty(\overline N),
\ea
has a bounded extension
\ba
D_{\tilde g}:H^1(N,|\tilde g|^{1/2}dx)\to \M(N;\R^3),
\ea
where $ \M(N;\R^3)$ denotes the space of $\quad\R^3$-valued signed Borel
measures on $N$.  Moreover, for $\tilde u\in X$,
we have,  in the sense of Borel measures
\beq\label{eq: sigma has zero measure}
(D_{\tilde g}\tilde u)(\Sigma)=0.
\eeq
\end{lemma}
\noindent

{\bf Proof.} Let $\tilde \phi\in  C^\infty(\overline N)$
and $\tilde \eta\in C(\overline N)$. Then $D^j_{\tilde g}\tilde
\phi\in L^\infty(N)$.
Let $\phi=\tilde \phi \circ F,\, \eta=\tilde \eta \circ F \in
L^\infty(\Omega)$.
Then,
\ba
\int_N (D^j_{\tilde g}\tilde \phi)\,\tilde \eta\,dx&=&
\int_{N\setminus \Sigma} (D^j_{\tilde g}\tilde \phi)\,\tilde \eta\,dx\\
&=&\int_{M_1\setminus \gamma_1}
|g|^{1/2}g^{kl}\frac{\p y^j}{\p x^l}\p_k\phi\,\eta\,dx+
\int_{M_2} |g|^{1/2}g^{kl}\frac{\p y^j}{\p x^l}\p_k\phi\,\eta\,dx.
\ea
As the  metric $g$ is bounded from above and below, and
$\frac{\p y^j}{\p x^l}=O(r^{-1})$ on $M_1$ and $=\delta^j_l$ on $M_2$,
we have
\ba
\left|\int_N (D^j_{\tilde g}\tilde \phi)\,\tilde \eta\,dx\right|
&\leq&
C_0(\|\phi\|_{H^1(M_1,dx)}\|\eta / r \|_{L^2(M_1,dx)}+
\|\phi\|_{H^1(M_2,dx)}\|\eta\|_{L^2(M_2,dx)})\\
&\leq&
C_1 ||\tilde \phi||_X ||\tilde \eta||_{C(N)}
  d^{1/2}_{\tilde g}(\hbox{supp}(\tilde \eta), \Sigma),
\ea
where $d_{\tilde g}$ is the distance on $N$ with respect to the metric
$\tilde g$.
This shows
the existence of the bounded extension
$D_{\tilde g}:H^1(N,|\tilde g|^{1/2}dx)\to \M(N;\R^3)$.
Also, if we consider functions $\tilde \eta$ supported
in small neighborhoods of $\Sigma$, we see that
(\ref{eq: sigma has zero measure}) follows.\nolinebreak
\hfill\proofbox

We also need  the following auxiliary result

\begin{lemma}\label{lem: aux}
Assume that
$\tilde u$
is a measurable function on $N$ such that
\beq\label{cond 1 A}
& &\tilde u\in L^2(N,\, |\tilde g|^{1/2}dx),\\
\label{cond 2 A}
& &\tilde u|_{N\setminus \Sigma}\in H^1_{loc}(N\setminus \Sigma,dx),
\\
\label{cond 3 A}
& &\int_{N\setminus\Sigma} |\tilde g|^{1/2}\tilde g^{ij}\p_i\tilde
u\p_j\tilde u\,dx<\infty.
\eeq
Then $\tilde u\in H^1(N,|\tilde g|^{1/2}dx)$.
\end{lemma}

Note that, due to the fact that $\tilde g$ is bounded and positive
definite on any
compact subset of $N\setminus \Sigma$,
condition (\ref{cond 2 A}) in fact follows from conditions (\ref{cond 1
A}),  (\ref{cond 3 A}) and is included for the convenience of future
references.

{\bf Proof.} Consider first the case when $\tilde u=0$ in
$N_1$.

First, the condition (\ref{cond 3 A})
implies that ${\tilde v}={\tilde u}|_{N_2}\in H^1(N_2,dx)$.
Let $f=v|_{\Sigma}\in H^{1/2}(\Sigma)$ and
$E^f\in H^1(N_1, dx)$ be an extension
of $f$. Let $\chi\in C^\infty_0(\R)$
be a cut-off function with $\chi(t)= 1$ for $|t|<\frac12$  and
$\chi(t)=0$ for
$|t|>1$.
We introduce Fermi coordinates near $\Sigma$ as in \S\ref{gbc},
$(\tau, \omega), \, \tau \in (0,2), \, \omega=(\omega_1, \omega_2) \in \Sigma$.

Define, for $\e>0$,
\ba
w_\e(x)=\left\{\begin{array}{cl}
v(x),& x\in N_2,\\
\chi(\frac {\tau}{\e}) E^f(x),& x\in N_1.
\end{array}\right.
\ea
Then $w_\e\in H^1(N,dx)$ and,
using (\ref{F1}), (\ref{propertyextra}), we see that
\beq
\label{additional}
\lim_{\e\to 0}\int_{N\setminus\Sigma} |\tilde g|^{1/2}
[\tilde g^{ij}\p_i(w_\e-\tilde u)\p_j
(w_\e-\tilde u)+(w_\e-\tilde u)^2]\,dx
\\ \nonumber
=\lim_{\e\to 0}\int_{N_1} |\tilde g|^{1/2}
[\tilde g^{ij}\p_iw_\e \p_j
w_\e+|w_\e|^2]\,dx
=0,
\eeq
Observe that the integrand vanishes
outside the a neighborhood of $\Sigma^+$ of
volume less than $C\e$.
Next, divide
the integral involving derivatives in the right-hand side of
(\ref{additional}) into the terms involving components tangential and
normal to the boundary,  using the fact that $\tau=2(r-1)$:
\ba
\int_{N_1\setminus \Sigma} |\tilde g|^{1/2} \chi^2(\frac {\tau}{\e})
{\tilde g}^{\alpha \beta} \p_{\omega_\alpha}E^f\,\p_{\omega_\beta}E^f\,
d\tau d\omega_1 d\omega_2,
\ea
and where $\alpha,\beta$ run over $\{1,2\}$,
\ba
\int_{N_1\setminus \Sigma} |\tilde g|^{1/2}
\left|\p_\tau[\chi(\frac {\tau}{\e}) E^f]\right|^2\,
d\tau d\omega_1 d\omega_2.
\ea
As, by (\ref{propertyextra}), $|\tilde g|^{1/2} {\tilde g}^{\alpha \beta}$
is bounded, the integral involving tangential derivatives tends to $0$
due to the volume of the domain of integration. Again, by
(\ref{propertyextra}) we have $|\tilde g|^{1/2}
\leq C\tau^2$; this, together with the volume estimate
and the fact that $|\p_\tau\chi(\frac {\tau}{\e})| \leq C \tau^{-1}$,
implies that the
integral  involving normal derivatives tends to $0$ when $\e \to 0$.
Similarly, we see that
\ba
\int_{N_1\setminus \Sigma} |\tilde g|^{1/2} |\chi(\frac {\tau}{\e})
E^f|^2 dx \to 0 \quad
\hbox{for}\,\, \e \to 0.
\ea
The function $w_\e\in H^1(N,dx)$ can be approximated with an arbitrarily
small error in  $H^1(N,dx)$ by a $C^\infty(\overline N)$
function,  and we see that the same holds in the $X$-norm.
Thus  $w_\e\in H^1(N,|\tilde g|^{1/2}dx)$, and
the above limit shows that
${\tilde u}\in   H^1(N,|\tilde g|^{1/2}dx)$.

Now let $\tilde u$ be a measurable function in $N$
satisfying (\ref{cond 1 A}), (\ref{cond 2 A}), and (\ref{cond 3 A}).
Let $\chi_{N_2}$ be the characteristic function of $N_2$.
As  $\chi_{N_2}\tilde u\in H^1(N,|\tilde g|^{1/2}dx)$,
it is enough to show that
            $\tilde u-\chi_{N_2}\tilde u\in H^1(N,|\tilde g|^{1/2}dx)$.
This means that it is enough to consider the case when
            $\tilde u=0$ in $N_2$. Clearly, we can restrict our attention
to the case when $\tilde u$
vanishes also near $\p N$.

Now let $u_1=\tilde u\circ F$ in $M_1\setminus \gamma_1 $. Then we see that
\ba
\int_{M_1\setminus\gamma_1} | g|^{1/2}g^{ij}\p_i(u_1)\p_j (u_1)
\,dx<\infty.
\ea

Let $w=\nabla u|_{M_1\setminus \gamma_1}$.
Using a change of coordinates in integration and (\ref{cond 1 A}),
we see that $u\in L^2(M_1\setminus \gamma_1,dx)$.
{Extending $u_1$ and $w$
to functions $u_1^e$ and $w^e$ on
$\gamma_1$, }
we obtain functions
$u_1^e\in L^2(M_1,dx)$ and $\R^3$--valued function $w^e\in L^2(M_1,dx)$.
Now $\nabla u_1^e-w^e\in H^{-1}(M_1,dx)$ is supported on $\gamma_1$.
Since there are no  non-zero $H^{-1}(M_1,dx)$ distributions supported
on $\gamma_1$,
we see that $\nabla u_1^e=w^e\in L^2(M_1,dx)$.
Thus we see that  $u_1^e\in H^1(M_1,dx)$. In the following we identify
$u_1$ and $u_1^e$. As $u_1$ vanishes near $\p M_1$, and $\gamma_1$
consists of a single point {and thus is a $(2,1)$-polar set
\cite[pp.393--7]{Ma}}, there are
$\phi_j\in C^\infty_0(M_1\setminus \gamma_1)$ such that
$\phi_j\to u_1$ in $H^1(M_1,dx)$ as $j\to \infty$, that is,
\ba
\lim_{j\to \infty}\int_{M_1}
|g|^{1/2}[g^{ik}\p_i(\phi_j- u)\p_k
(\phi_j-u)+(\phi_j-u)^2]\,dx=0.
\ea
Now let $\tilde \phi_j\in C^\infty_0(N)$, with
$\hbox{supp}(\tilde \phi_j) \subset N_1$ and
\bfo
\tilde \phi_j=
\left\{\begin{array}{cl}
\phi_j\circ F_1^{-1}& \hbox{in } N_1,\\
0& \hbox{in } N_2.
\end{array}\right.
\efo
Then the previous equation implies that
\ba
& &\lim_{j\to \infty}\int_{N\setminus \Sigma}|\tilde g|^{1/2}
[\tilde g^{ik}\p_i(\tilde\phi_j-\tilde u)\p_k
(\tilde\phi_j-\tilde u)+(\tilde\phi_j-\tilde u)^2]\,dx =\\
\nonumber
& & \lim_{j\to \infty}\int_{N_1}|\tilde g|^{1/2}
[\tilde g^{ik}\p_i(\tilde\phi_j-\tilde u)\p_k
(\tilde\phi_j-\tilde u)+(\tilde\phi_j-\tilde u)^2]\,dx
=0,
\ea
where we use that $\tilde u = 0$ in $N_2$.

This shows that $\tilde \phi_j$ is a sequence converging in the $X$-norm
and that the limit is $\tilde u$.
Thus  $\tilde u\in   H^1(N,|\tilde g|^{1/2}dx)$,
proving the claim.\hfill\proofbox.
\medskip

{Although in this section  $(M,N,F,\gamma,\Sigma,g)$ continues to
denote a  single coating, we will see later that the following definition is
also appropriate for the double coating construction.}

Let $\tilde f\in L^2(N,dx)$ be a function such that $\supp (\tilde
f)\cap \Sigma=\emptyset$.

\begin{definition}
\label{energysolution}
{Let $(M,N,F,\gamma,\Sigma,g)$ be
a coating construction.}
A measurable function $\tilde u$ on $N$ is a
\emph{ finite energy }  solution of the Dirichlet problem for the Helmholtz
equation on
$N$,
\beq\label{case 3 B}
& &(\Delta_{\tilde g}+k^2)\tilde u=\tilde f\quad\hbox{on }N,\\
& &\tilde u|_{\p N}=\tilde h,\nonumber
\eeq
if
\beq\label{cond 1 B}
& &\tilde u\in L^2(N,\, |\tilde g|^{1/2}dx);\\
\label{cond 2 B}
& &\tilde u|_{N\setminus \Sigma}\in H^1_{loc}(N\setminus \Sigma,dx);
\\
\label{cond 3 B}
& &\int_{N\setminus\Sigma} |\tilde g|^{1/2}\tilde g^{ij}\p_i\tilde
u\p_j\tilde u\,dx<\infty,\\
& &\tilde u|_{\p N}=\tilde h,\nonumber;
\eeq
and, for all ${\tilde \psi} \in C^\infty(N)\hbox{ with } {\tilde
\psi}|_{\p N}=0$,
\beq\label{cond 5 B}
& &\int_{N} [-(D_{\tilde g}\tilde u) \p_j {\tilde \psi}+k^2 \tilde
u {\tilde \psi}|\tilde g|^{1/2}]dx
=\int_{N} \tilde f(x) {\tilde \psi} (x)
|\tilde g|^{1/2} dx
\eeq
where the integral on the left hand side of (\ref{cond 5 B})
is defined by distribution-test function duality.
\end{definition}

Note as before that condition (\ref{cond 2 B}) follows from (\ref{cond 1 B}),
(\ref{cond 3 B}).

Cloaking by the single coating of arbitrary active devices, with
respect to solutions of  the
Helmholtz equation at all frequencies, then follows from
the following.
\bigskip

\begin{theorem}\label{single coating with Helmholtz}
Let $u=(u_1, u_2) : M\setminus \gamma\to \R$
and  $\tilde u: N \setminus \Sigma\to \R$ be  measurable functions
such that $u=\tilde u\circ F$.
Let $f=(f_1, f_2): M \setminus \gamma\to \R$
and  $\tilde f:N\setminus \Sigma\to \R$ be $L^2$ functions
supported away from $\gamma$ and $\Sigma$
such that  $f=\tilde f\circ F$,
and $\tilde h:\p N\to\R,\, h:\p M_1\to\R$ be such that $h=\tilde
h\circ F_1$.

Then the following are equivalent:
\begin{enumerate}

\item The function
$\tilde u$, considered as a measurable function on $N$,
is a finite energy solution to the Helmholtz equation (\ref{case  3 B})
with inhomogeneity $\tilde f$ and Dirichlet data $\tilde  h$ in the sense of
Definition
\ref{energysolution}.

\item The function $u$ satisfies
\beq\label{eq on M1}
(\Delta_{g}+k^2)u_1=f_1\quad\hbox{ on }M_1,\quad
u_1|_{\p M_1 }=h,
\eeq
and
\beq\label{eq on M2}
& &(\Delta_{g}+k^2)u_2=f_2\quad\hbox{ on } M_2, \quad
g^{jk}\nu_j\p_k u_2|_{\p M_2}=b,
\eeq
with $b=0$. Here $u_1$ denotes the continuous extension of $u_1$ from
$M_1 \setminus \gamma$ to $M_1$
\end{enumerate}

Moreover, if $u$ solves (\ref{eq on M1}) and  (\ref{eq on M2})  with
$b\not =0$, then the function
$\tilde u=u\circ F^{-1}:N\setminus \Sigma\to \R$,
considered as a measurable function on $N$,
is not a finite energy solution to the  Helmholtz equation.
\end{theorem}

{\bf Remarks.} (i) It follows that the construction of
\cite{GLU1,PSS1} cloaks active devices
from detection by unpolarized EM waves at all frequencies.

(ii) Observe that in (\ref{eq on M1}) the right
hand side $f_1$ is zero near
         $\gamma_1$.
Thus $u_1$, considered as a distribution in a neighborhood of $\gamma_1$,
has an extension on $\gamma_1$ that is  $C^\infty$ smooth function in
a neighborhood of  $\gamma_1$.

(iii) As noted previously, for
the single coating case
one may assume that $N_2=M_2$ and
$F|_{M_2}$ is the identity. Thus $\tilde u|_{N_2}=u|_{M_2}$; hence, if
$\tilde u$
is a finite energy solution of  the Helmholtz equation on $N$, we
see that  $u|_{M_2}$ satisfies the Neumann boundary condition on  $\p
M_2$ and thus
        also $\tilde u|_{N_2}$ automatically has to satisfy the Neumann
condition on $\Sigma^-$.
The Neumann boundary condition
that appears  on $\p N_2$ means that, observed from the inside of the
cloaked region $N_2$,
the single coating construction  has the effect of creating  a virtual
sound hard, i.e., perfectly reflecting, surface at $\Sigma$.
Similarly, we will see later that there are hidden
boundary conditions for Maxwell's equations in the presence of the single
coating, but they are overdetermined and generally preclude such solutions
existing.

\noindent
{\bf Proof.} First we proof that Helmholtz on $M$
implies Helmholtz on $N$.

Let $f\in L^2(M,dx)$ be a function such that
$\supp (f)\cap (\gamma
\cup \p M_1\cup \p M_2)=\emptyset$.
Assume that a function $u$ on $M$ is a classical solution of (\ref{eq on
M1}) and (\ref{eq on M2}).
Notice that we have required here
that $u_2$ on $\p M_2$ satisfies the Neumann boundary condition at $\p M_2$.

Again, define $\tilde u=F_*u$  and $\tilde f=f\circ F^{-1}$ on
$N\setminus \Sigma$ and extend it, e.g., by  setting it equal to zero on
$\Sigma$. Note that then
$\tilde f\in L^2(N,dx)$ is supported away from
$\Sigma$, and
$\tilde u\in  L^2(N,|\tilde g|^{1/2}dx)$
            satisfies
\beq\label{case 3B 2}
(\Delta_{\tilde g}+k^2)\tilde u_1=\tilde f_1= \tilde f|_{N_1}\quad\hbox{in
}N_1,\quad
\tilde u|_{\p N}=\tilde h,
\eeq
and
\beq\label{case 3C 2}
& &(\Delta_{\tilde g}+k^2)\tilde u_2=\tilde f_2= \tilde
f|_{N_2}\quad\hbox{in }N_2.
\eeq

Let $\Sigma(\e)$ be the $\e$-neighborhood
of $\Sigma$  with respect to the metric  $\tilde g$.
Let $\gamma(\e)$ be the $\e$-neighborhood
of $\gamma\subset M_1$ with respect to the metric $g$.
Let $g_{bnd}$ and  $\tilde g_{bnd}$ be the induced metrics
on $\p \gamma(\e)$ and $\p \Sigma(\e)$, correspondingly.

Clearly, the function $\tilde u$ satisfies conditions (\ref{cond 1 B}),
(\ref{cond 2 B}), and  (\ref{cond 3 B}).
By Lemma \ref{lem: aux}, we have that
$\tilde u\in H^1(N,|\tilde g|^{1/2}dx)$, and
$D_{\tilde g}\tilde u$ is thus well defined.

Using relations (\ref{propertyextra}) for the normal
component and (\ref{case 3B 2}), (\ref{case 3C 2}),
{and property (\ref{eq: sigma has zero measure})
of $D_{\tilde g} u$,}
we see that, for $\tilde \psi\in C^\infty_0(N)$,
\beq\label{eq: weak N 2 B}
& &\int_{N} [-D_{\tilde g}(\tilde u) \p_j \tilde
\psi+k^2 \tilde u\tilde \psi|\tilde g|^{1/2}-\tilde f\tilde
\psi|\tilde g|^{1/2}]dx\\ \nonumber
&=&\lim_{\e\to 0} \int_{N\setminus \Sigma(\e)} (-\tilde
g^{ij}\,\p_i \tilde u\, \p_j \tilde \psi+(k^2 \tilde u-\tilde
f)\tilde \psi)|\tilde g|^{1/2}dx\\
            \nonumber
&=&\lim_{\e\to 0} (\int_{\p \Sigma(\e)\cap N_2}+\int_{\p
\Sigma(\e)\cap N_1}
)
            (-\tilde g^{ij}\,\nu _j\,\p_i \tilde u\,  \tilde \psi)|\tilde
g_{bnd}|^{1/2}dS
            \nonumber\\
&=&
\lim_{\e\to 0} \int_{\Sigma(\e) \cap N_2}(-\tilde g^{ij}\,\nu _j\p_i
\tilde u_2 \tilde \psi)
|\tilde g_{bnd}|^{1/2}dS+
            \label{int 1}\\
& &+
\lim_{\e\to 0} \int_{\p \gamma(\e)}
            (-g^{ij}\,\p_i u_1\, \nu _j (\tilde \psi\circ F))|g_{bnd}|^{1/2}dS
            \label{int 2}\\
&=&0 \nonumber.
\eeq

Indeed, the integral (\ref{int 1}) in the right-hand side of
this equation tends to $0$ due  to the boundary condition on
$\Sigma^-$  (\ref{eq on M2}),  and boundedness of
$\tilde \psi\circ F$. To analyze the integral (\ref{int 2}) observe that, as
$\hbox{supp} f_1 \cap \gamma_1 = \emptyset$, $u_1$ is infinitely smooth
near $\gamma_1$. Thus all $\p_iu_1$  and $\tilde \psi\circ F$ are
bounded near $\gamma_1$, while
            the area of
${\p \gamma(\e)}$ is bounded by $C\e^2$.
Hence we see that (\ref{cond 5 B}) is valid
and thus
\ba
(\Delta_{\tilde g}+k^2)\tilde u=\tilde f\quad\hbox{in }N
\ea
in the sense of the Definition \ref{energysolution}.

Summarizing, so far we have proven that a (classical) solution to the Helmholtz
equation on $M$ yields, via the pushforward, a finite energy solution
to the equation on $N$.

Next we consider the other direction and prove that
the Helmholtz equation on $N$ implies Helmholtz equation on $M$.

Assume that $\tilde u$ satisfies Helmholtz equation (\ref{case 3 B})
on $(N,\tilde g)$ in the sense of Definition \ref{energysolution}, with $\tilde
f\in L^2(N)$
supported away from $\Sigma$. In particular, $\tilde u$
is  a measurable function in $N$
satisfying (\ref{cond 1 A}),
(\ref{cond 2 A}),
and (\ref{cond 3 A}).

Let  $u=\tilde u\circ F$ and $f=\tilde f\circ F$ on $M\setminus \gamma$.
Then we have
\beq \label{eq on M_1 minus gamma}
(\Delta_g+k^2)u_1=f_1=f|_{M_1 \setminus \gamma_1} \quad \hbox{in }
M_1\setminus \gamma_1,\quad
u_1|_{\p M_1}=h
\eeq
and
\beq \label{eq on M_2}
(\Delta_g+k^2)u_2=f_2=f|_{M_2} \quad \hbox{in } M_2.
\eeq

By conditions (\ref{cond 1 A}),
(\ref{cond 2 A}),
and (\ref{cond 3 A}),
we have that
\ba
& & |u|^2\in L^1(M_1\setminus \gamma_1,|g|^{1/2}dx),\\
& & g^{jk}(\p_j u)(\p_k u)\in L^1(M_1\setminus \gamma_1,|g|^{1/2}dx).
\ea
and thus $u_1\in H^1(M_1\setminus \gamma_1,dx)$.
As before, we see that
\beq \label{eq on M_1}
(\Delta_g+k^2)u_1=f_1\quad \hbox{in } M_1,\quad
u_1|_{\p M_1}=h,
\eeq
where $f_1$ is extends to have the value $0$ at $\gamma_1$ and $u_1$ is
smooth near $\gamma_1$.

Let us now consider the boundary condition on $M_2$.
Since $\tilde u$ satisfies (\ref{cond 5 B}),
we  see that for ${\tilde \psi} \in C^\infty_0(N)$,
\beq\label{eq: weak N 2 c}
0&=&\int_{N} [-D_{\tilde g}\tilde u \p_j \tilde
\psi+k^2 \tilde u\tilde \psi |\tilde g|^{1/2}-\tilde f\tilde
\psi |\tilde g)|^{1/2}]dx
\\ \nonumber
&=&\lim_{\e\to 0} \int_{N\setminus \Sigma(\e)} (-\tilde
g^{ij}\,\p_i \tilde u\, \p_j \tilde \psi+(k^2 \tilde u-\tilde
f)\tilde \psi)\ |\tilde g|^{1/2}dx\\
\nonumber
&=&\lim_{\e\to 0} (\int_{\p \Sigma(\e)\cap N_2}+\int_{\p
\Sigma(\e)\cap N_1}
)
(-\tilde g^{ij}\,\nu _j\,\p_i \tilde u\,  \tilde \psi) |\tilde
g_{bnd}|^{1/2}dS
\nonumber\\
&=&
\int_{\p M_2}(-g^{ij}\,\nu _j\p_i u_2|_{\p M_2}   \psi)
|g_{bnd}|^{1/2}dS
\\ \nonumber
& &+\lim_{\e\to 0} \int_{\p \gamma(\e)}
(-g^{ij}\,\p_i u_1\, \nu _j  \psi)|g_{bnd}|^{1/2}dS\\
&=&  \nonumber
\int_{\p M_2}(-g^{ij}\,\nu _j\p_i u|_{\p M_2}   \psi)
|g_{bnd}|^{1/2}dS,
\eeq
where $\psi =\tilde \psi \circ F$. Here
we use  the fact that
$u_1$ is a smooth function, implying that $\p_i u_1$ is bounded
and that $\psi=\tilde \psi\circ F $ is  bounded.
As ${\tilde \psi}|_{\p M_2} \in C^{\infty}(\p M_2)$
is arbitrary,
this shows that
\beq\label{Neumann boundary condition on M_2}
\tilde g^{ij}\,\nu _j\p_i \tilde u_2|_{\p M_2} =0.
\eeq

Thus, we have shown that
the  function $u$ is a classical solution on $M$ of
\beq\label{eq on M1 B}
(\Delta_{g}+k^2)u_1=f_1 \quad\hbox{in }M_1,\quad
u_1|_{\p M_1}=h
\eeq
and
\beq\label{eq on M2 B}
(\Delta_{g}+k^2)u_2=f_2 \quad\hbox{in }M_2,\quad
g^{jk}\nu_j\p_k u_2|_{\p M_2}=0.
\eeq
This proves the claim, and finishes the proof of
Theorem \ref{single coating with Helmholtz}. \hfill\proofbox
\bigskip

\subsubsection{Operator theoretic definition of the Helmholtz equation}

It is standard in quantum physics that a self-adjoint operator
can be defined via the quadratic form corresponding to
energy. In the case  considered here, the energy associated with the wave
operator is defined by the quadratic (Dirichlet) form $A$,
\beq
\label{Dirichletform}
A[ \tilde u,  \tilde u]:=\int_{N\setminus\Sigma} |\tilde
g|^{1/2}\tilde g^{ij}\p_i\tilde
u\p_j\tilde u\,dx, \quad \tilde u\in {\mathcal D}(A)
\eeq

As we deal with the sound-soft boundary $\p N$ or,
more generally, with the source on $\p N$ of the form
$\tilde u|_{\p N}= \tilde h$, the domain ${\mathcal D}(A)$
of the form $A$ should be taken as
\ba
{\mathcal D}(A)=H^1_0(N, |\tilde g|^{1/2} dx)\subset X.
\ea
Thus, by  standard techniques of  operator theory, e.g., \cite{Kato},
the form $A$ defines a positive selfadjoint operator, denoted
$A_0=-\Delta^D_{\tilde g}$, on
$L^2(N, |\tilde g|^{1/2} dx)$. Next we recall this construction.
We say that
$\tilde u
\in H^1_0(N, |\tilde g|^{1/2} dx)$ is in
the domain of $A_0$, $\tilde u\in {\mathcal D}(A_0)$ if
there is an $\tilde f \in L^2(N,
|\tilde g|^{1/2} dx)$
such that for all $\tilde v \in
H^1_0(N, |\tilde g|^{1/2} dx)$,
\beq
\label{pre-energyidentity}
A[\tilde u,\tilde v]=
\int_{N} {\tilde f} \tilde v\,|\tilde g|^{1/2}  dx.
\eeq
In this case, we define
\ba
A_0\tilde u={\tilde f}.
\ea

\begin{proposition}
Assume that $- k^2$ is not in the spectrum of $\Delta^D_{\tilde g}$.
Then $\tilde u$ is {a finite energy} solution to
\bfo
(\Delta_{\tilde g}+k^2)\tilde u=\tilde f, \quad \tilde u|_{\p N}=
\tilde h \in H^{1/2}(\p N)
\efo
if and only if
\beq\label{new def.}
\tilde u=E\tilde h +(\Delta^D_{\tilde g}+k^2)^{-1}(\tilde f-
(\Delta_{\tilde g}+k^2)E\tilde h),
\eeq
where $E\tilde h$ is an $H^1(N,dx)$-extension
of $\tilde h$ to $N$  satisfying
$\hbox{supp}\,(E\tilde h) \subset \p N\cup N_1$.
\end{proposition}

\noindent
{\bf Proof.}
First we show that if $\tilde u$ satisfies the conditions of Definition
\ref{energysolution} then it satisfies (\ref{new def.}).
As ${\tilde \psi} \in C^\infty(N),\, {\tilde \psi}|_{\p N}=0$,
imply that  ${\tilde \psi} \in H^1_0(N, |g|^{1/2} dx)$,
we see  by (\ref{cond 5 B})  that $ \tilde u -E\tilde h$ satisfies
\ba
      \int_{N} (-D^j_{\tilde g}(\tilde u -E\tilde h)\, \p_j\tilde v
+k^2(\tilde u -E\tilde h) \tilde v)\,dx
=
\int_{N\setminus\Sigma} |\tilde g|^{1/2} ({\tilde f}- (\Delta_{\tilde
g}+k^2)E\tilde h)
\tilde v\,dx,
\ea
for any  $\tilde v \in C^{\infty}_0(N)$.
By (\ref{eq: sigma has zero measure}) and  (\ref{cond 3 B}),
this implies
\beq\label{aim 1}
& &  \int_{N\setminus\Sigma} |\tilde g|^{1/2}\left(-\tilde
g^{ij}\p_i(\tilde u -E\tilde h)\p_j\tilde v
+k^2(\tilde u -E\tilde h) \tilde v\right)\,dx
\\ \nonumber
& &=
\int_{N\setminus\Sigma} |\tilde g|^{1/2} ({\tilde f}- (\Delta_{\tilde
g}+k^2)E\tilde h)
\tilde v\,dx,
\eeq
for any $\tilde v \in C^\infty_0(N)$.
We need to show that (\ref{aim 1}) is valid for all
$\tilde v \in H^1_0(N, |g|^{1/2} dx)$.

Observe that
\bfo
\int_{N\setminus\Sigma} |\tilde g|^{1/2}\left(-\tilde
g^{ij}\p_i(E\tilde h)\p_j\tilde v +k^2(E\tilde h) \tilde v\right)\,dx =
\int_{N\setminus\Sigma} |\tilde g|^{1/2} ((\Delta_{\tilde g}+k^2)E\tilde h)
\tilde v\,dx,
\efo
where we use that $\hbox{supp}(E\tilde h) \subset \p N \cup N_1$ and $\tilde
v|_{\p N}=0$.
Thus, it remains to show that
\beq\label{eq: intermediate step}
          \int_{N\setminus\Sigma} |\tilde g|^{1/2}\left(-\tilde
g^{ij}\p_i \tilde u \p_j\tilde v
+k^2 \tilde u  \tilde v\right)\,dx =
\int_{N\setminus\Sigma} |\tilde g|^{1/2} \tilde f
\tilde v\,dx
\eeq
for $\tilde v \in H^1_0(N, |g|^{1/2} dx)$.
Clearly, to show this
it is enough to show that
\beq\label{eq: lim}
\lim_{\e\to 0}  \int_{N\setminus\Sigma_1(\e)} |\tilde g|^{1/2}\left(-\tilde
g^{ij}\p_i \tilde u \p_j\tilde v
+k^2 \tilde u  \tilde v- \tilde f \tilde v\right)\,dx =0.
\eeq
where $\Sigma_1(\e)=N_1\cap\Sigma(\e)$.

Next we argue analogously to the reasoning that led
to equation (\ref{eq: weak N 2 B}).
Let $v=\tilde v\circ F$, $f=\tilde f\circ F$, and
$u=\tilde u\circ F$ in $M\setminus \gamma$.
To clarify notations, denote $u_1=u|_{M_1}$, $u_2=u|_{M_2}$
, $v_1=v|_{M_1}$, $v_2=v|_{M_2}$,
and  $f_1=f|_{M_1}$, $f_2=f|_{M_2}$.
Then, by Proposition \ref{single coating with Helmholtz},
\beq\label{recall 1}
& &(\Delta_g+k^2)u_1=f_1,\quad \hbox{in }M_1,\\
& &(\Delta_g+k^2)u_2=f_2,\quad \hbox{in }M_2, \label{recall 2}\\
& &\p_\nu u_2|_{\p M_2}=0, \label{recall 3}
\eeq
and we see that
\ba
& & \lim_{\e\to 0} \int_{N\setminus\Sigma_1(\e)} |\tilde g|^{1/2}\left(-\tilde
g^{ij}\p_i \tilde u \p_j\tilde v
+k^2 \tilde u  \tilde v- \tilde f \tilde v\right)\,dx \\
&=&  \lim_{\e\to 0} \int_{(M_1\setminus\gamma(\e))\cup M_2} |g|^{1/2}\left(-
g^{ij}\p_iu \p_j v
+k^2 uv-  f  v\right)\,dx
\\
&=&  \int_{\p \gamma(\e)}
         (-g^{ij}\nu _j\,\p_i
         (u)\,  v)|
g_{bnd}|^{1/2}dS+   \int_{\p  M_2}
         (-g^{ij}\nu _j\,\p_i
         (u)\,  v)|
g_{bnd}|^{1/2}dS.
\ea
By (\ref{recall 3}), we have that
\beq\label{eq: lim 2}
       \int_{\p M_2}
         (-g^{ij}\nu _j\,\p_i
         (u)\,  v)|\tilde
g_{bnd}|^{1/2}dS=0.
\eeq
Next we consider
\ba
I_1(\e)=\int_{\p \gamma(\e)\cap M_1}
         (-g^{ij}\nu _j\,\p_i
(u)\,  v)|\tilde
g_{bnd}|^{1/2}dS.
\ea
Note that   $\lim_{\e  \to 0} I_1(\e)$ exists
as the limits (\ref{eq: lim}) and (\ref{eq: lim 2})
exists.

As $\supp(f)\cap \gamma=\emptyset $,  we see that $u_1$ is smooth function
near $\gamma$.
Moreover, as $\tilde v\in X$,
we observe that $v_1 \in H^1(M_1\setminus \gamma, dx)$,
and so $v_1$ can be extended to   $v_1 \in H^1(M_1, dx)$.
Hence, by the Sobolev embedding
theorem,
$v_1 \in L^6(M_1, dx)$. This allows us to deduce that
\beq
\label{l-6}
\liminf_{\e \to 0} \e^{-3/2} \int_{\p \gamma(\e)} |v_1|\quad dS =0.
\eeq
Indeed,
\ba
& &
\int_0^\e \left(\int_{\p  \gamma(r)} |v_1| dS(x)\right) dr =
\int_{ \gamma(\e)}
|v_1|dx
\\ \nonumber
& & \leq
\left(\int_{ \gamma(\e)} |v_1|^6dx\right)^{1/6}\,\left(\int_{ \gamma(\e)}
dx\right)^{5/6}= o(\e^{5/2}).
\ea
Clearly, this inequality implies (\ref{l-6}). Thus using
boundedness of $g^{ij}\nu_j\p_i u$ we see that
\ba
\liminf_{\e  \to 0} \int_{\p \gamma(\e)}
            (- g^{ij}\nu _j\,\p_i ( u)\,   v |
g_{bnd}|^{1/2})dS =0.
\ea
As  $\lim_{\e  \to 0} I_1(\e)$ exists, this implies
$\lim_{\e  \to 0} I_1(\e)=0$.
As $\tilde u|_{\p N}= \tilde h$ by
Definition \ref{energysolution}
we have shown that
Definition \ref{energysolution} implies (\ref{new def.}).

Next, consider the case when $\tilde u$ satisfies (\ref{new def.}).
Since  $\tilde u\in X$,
we see by  (\ref{eq: sigma has zero measure}) that
\beq\label{eq: pre-intermediate}
          \int_{N} D_{\tilde g}^j(\tilde u) \p_j\tilde v \,dx
=  \int_{N\setminus\Sigma} |\tilde g|^{1/2}\tilde
g^{ij}\p_i \tilde u \p_j\tilde v \,dx
\eeq
for all $\tilde v\in C^\infty_0(N)$.
Thus, by (\ref{new def.}) we see that (\ref{eq: intermediate step})
is valid for
         $\tilde v \in C^\infty_0(N)$, which implies
condition (\ref{cond 3 B}). The other
conditions in
Definition \ref{energysolution} follow easily from (\ref{new def.}).
\hfill\proofbox

\subsection{Helmholtz for the double coating}

{We now examine solutions to the Helmholtz equation in the presence of
the double
coating; we will establish
full-wave invisibility at all nonzero frequencies. Unlike for the
single coating,
for the double coating  no extra boundary conditions appear at  $\Sigma$.
Otherwise, the reasoning here parallels that in \S 3.2.}

Throughout this section,  $(M,N,F,\gamma,\Sigma,g)$ is
a double coating construction.

\subsubsection{Weak solutions for the double coating.}

Suppose that $k\geq 0$ and $\tilde f\in L^2 (N, |\tilde g|^{1/2}dx)$.
We use the same notion of weak solution as for the single coating,
saying that $\tilde u$ is a
\emph{finite energy} solution of
\beq\label{case 3}
(\Delta_{\tilde g}+k^2)\tilde u=\tilde f\quad\hbox{in }N,\quad
\tilde u|_{\p N}=\tilde h\eeq
if $\tilde u$
is a solution of the Dirichlet problem in the sense  of
Definition \ref{energysolution}.

We start with analogues of the space $H^1(N, |\tilde g|^{1/2}dx)$, and
Lemmas \ref{lem: extension} and \ref{lem: aux}.
To this end define, for ${\tilde \phi}\in  C^\infty(\overline N)$,
\ba
\|{\tilde \phi} \|_Y^2:=\int_{N}
(|\tilde g|^{1/2}\tilde g^{ij}\p_i{\tilde \phi} \p_j{\tilde \phi}+
|\tilde g|^{1/2}|{\tilde \phi}|^2)\,dx.
\ea
Let
\ba
H^1(N,|\tilde g|^{1/2}dx)=Y:=\hbox{cl}_{Y}( C^\infty(\overline N))
\ea
be the completion of  $ C^\infty(\overline N)$ with respect to the norm
$\|\,\cdotp\|_Y$. Note that $H^1(N,|\tilde g|^{1/2}dx)\subset
L^2(N,|\tilde g|^{1/2}dx)$, so we can consider its elements
as measurable functions in $N$.

\begin{lemma}\label{lem:3.7a}
The map
\ba
\phi\longrightarrow D_{\tilde g}{\tilde \phi}=(D^j_{\tilde g}{\tilde
\phi})_{j=1}^3=
(|\tilde g|^{1/2}\tilde g^{ij}\p_i{\tilde \phi})_{j=1}^3,\quad \phi\in
C^\infty(\overline N),
\ea
has a bounded extension
\ba
D_{\tilde g}:H^1(N,|\tilde g|^{1/2}dx)\to \M(N;\R^3),
\ea
where $ \M(N;\R^3)$ denotes the space of $\quad\R^3$-valued signed Borel
measures on $N$.  Moreover, for $\tilde u\in Y$,
we have
\beq\label{analog12}
(D_{\tilde g}\tilde u)(\Sigma)=0.
\eeq
If $\tilde u$
is a measurable function on $N$ such that
\beq\label{L3.71}
& &\tilde u\in L^2(N,\, |\tilde g|^{1/2}dx),\\
\label{L3.72}
& &\tilde u|_{N\setminus \Sigma}\in H^1_{loc}(N\setminus \Sigma,dx),
{\hbox{ and }}
\\
\label{L3.73}
& &\int_{N\setminus\Sigma} |\tilde g|^{1/2}\tilde g^{ij}\p_i\tilde
u\p_j\tilde u\,dx<\infty,
\eeq
then $\tilde u\in H^1(N,|\tilde g|^{1/2}dx)$.
\end{lemma}

\noindent {\bf Proof.}
The proof here is essentially the same as of
Lemmas \ref{lem: extension} and \ref{lem: aux}.
The only difference is that, as described in \S\ref{gbc}, the map
\ba
F: M\setminus \gamma \rightarrow N\setminus\Sigma
\ea
now consists of two maps,
\ba
F_i:  M_i\setminus \gamma\rightarrow N_i ,
\quad i=1,2,
\ea
{having similar structure to each other, namely  that of the
map $F_1$ in
the single coating construction. (Recall that for the double coating
construction,
$\gamma_1:=\gamma\cap M_1$ is  a point $O \in M_1$ and
$\gamma_2:=\gamma\cap M_2$  a point $NP \in M_2$.)}

Therefore, when proving  that $\tilde u$ satisfying
(\ref{L3.71})--(\ref{L3.73}) is in $H^1(N,|\tilde g|^{1/2}dx)$,
we can use the fact that, in this case, both
$(1-\chi_{N_1})\tilde u$ and $(1-\chi_{N_2})\tilde u$
satisfy (\ref{L3.71})--(\ref{L3.73}) and carry out the proof
for each of them as for the  $(1-\chi_{N_2})\tilde u$ term
in the proof of Lemma \ref{lem: aux}.

Invisibility of active devices in the
presence of the double coating with respect to the Helmholtz equation
at all frequencies
then follows from

\begin{theorem}\label{double coating with Helmholtz:a}
Let $u=(u_1, u_2):M\setminus \gamma\to \R$
and  $\tilde u:N\setminus \Sigma\to \R$ be  measurable functions
such that $u=\tilde u\circ F$.
Let $f=(f_1, f_2):M\setminus \gamma\to \R$
and  $\tilde f:N\setminus \Sigma\to \R$ be $L^2$ functions
supported away from  $\gamma$ and $\Sigma$
such that  $f=\tilde f\circ F$.
Then the following are equivalent:

\begin{enumerate}
\item The function
$\tilde u$, considered as a measurable function on $N$,
is a finite energy solution to the Helmholtz equation (\ref{case 3}) with
inhomogeneity $\tilde f$ and Dirichlet data $\tilde h$ in the sense of
Definition
\ref{energysolution}.
\item We have
\beq\label{analog50}
(\Delta_{g}+k^2)u_1=f_1\quad\hbox{ on }M_1,\quad
u|_{\p M}=h:= \tilde h \circ F
\eeq
and
\beq\label{analog51}
& &(\Delta_{g}+k^2)u_2=f_2\quad\hbox{ on }M_2.
\eeq
\end{enumerate}
\end{theorem}

\noindent {\bf Proof}
As in the proof of Theorem \ref{single coating with Helmholtz},
we first prove that Helmholtz on $M$ implies Helmholtz on $N$.

Let $f\in L^2(M,dx)$ be a function such that
$\supp (f)\cap (\gamma
\cup \p M_1\cup \p M_2)=\emptyset$.
Assume that a function $u=(u_1, u_2)$ on $M$ is a
classical solution of (\ref{analog50})
and (\ref{analog51}).
Define $\tilde u=F_*u$  and $\tilde f=f\circ F^{-1}$ on
$N\setminus \Sigma$ and extend it, e.g., by  setting it equal to
zero on $\Sigma$.
Note that then $\tilde f\in L^2(N,dx)$ is supported away of $\Sigma$.
Then
$\tilde u\in
L^2(N,|\tilde g|^{1/2}dx)$
satisfies

\beq\label{analog24}
(\Delta_{\tilde g}+k^2)\tilde u_1=\tilde f_1= \tilde f|_{N_1}\quad\hbox{in
}N_1,\quad
\tilde u|_{\p N}=\tilde h,
\eeq
and
\beq\label{analog25}
& &(\Delta_{\tilde g}+k^2)\tilde
u_2=\tilde f_2= \tilde
f|_{N_2}\quad\hbox{in }N_2.
\eeq

Let
$\Sigma(\e)$ be the $\e$-neighborhood
of $\Sigma$  with respect to
the metric  $\tilde g$.
Let $\gamma_1(\e)$ be the
$\e$-neighborhood
of $\gamma_1=\{0\}\subset M_1$ with respect to the
metric $g$.
Let $\gamma_2(\e)$ be the $\e$-neighborhood
of
$\gamma_2=\{NP\}\subset M_2$ with respect to the metric $g$.
Let
$g_{bnd}$ and  $\tilde g_{bnd}$ be the induced metrics
on $\p
\gamma(\e)$ and $\p \Sigma(\e)$, correspondingly.

Clearly, the
function $\tilde u$ satisfies conditions (\ref{cond 1 B}),
(\ref{cond
2 B}), and  (\ref{cond 3 B}).
By Lemma \ref{lem:3.7a}, we have
that
$\tilde u\in H^1(N,|\tilde g|^{1/2}dx)$, and
$D_{\tilde g}\tilde
u$ is thus well defined.

Using relations (\ref{propertyextra}),
(\ref{property 1})
in $M_1$ and (\ref{property 1 b}) in
$M_2$, it
follows from
(\ref{analog24}), (\ref{analog25}) that
for $\tilde
\psi\in C^\infty_0(N)$,
\beq\label{analog26}
& &\int_{N} [-(D_{\tilde
g})\tilde u \p_j \tilde
\psi+k^2 \tilde u\tilde \psi|\tilde
g|^{1/2}-\tilde f\tilde
\psi|\tilde g|^{1/2}]dx\\
\nonumber
&=&\lim_{\e\to 0} \int_{N\setminus \Sigma(\e)}
(-\tilde
g^{ij}\,\p_i \tilde u\, \p_j \tilde \psi+(k^2 \tilde
u+\tilde
f)\tilde \psi)|\tilde g|^{1/2}dx\\
\nonumber
&=&\lim_{\e\to
0} (\int_{\p \Sigma(\e)\cap N_2}+\int_{\p
\Sigma(\e)\cap
N_1}
)
(-\tilde g^{ij}\,\nu _j\,\p_i \tilde u\,  \tilde
\psi)|\tilde
g_{bnd}|^{1/2}dS
\nonumber\\
&=&
\lim_{\e\to 0} \int_{\p
\gamma_1(\e)}
(-g^{ij}\,\p_i u_1\, \nu _j (\tilde \psi\circ
F))|g_{bnd}|^{1/2}dS
\nonumber\\
& &+
\lim_{\e\to 0} \int_{\p
\gamma_2(\e)}
(-g^{ij}\,\p_i u_2\, \nu _j (\tilde \psi\circ
F))|g_{bnd}|^{1/2}dS
\nonumber\\
&=&0 \nonumber.
\eeq
Indeed, both
terms in the right-hand side of (\ref{analog26})
tend to $0$ by the
same arguments as the term
$\int_{\p \gamma(\e)}
(-g^{ij}\,\nu
_j\,\p_i u_1\,  (\tilde \psi\circ F))|g_{bnd}|^{1/2}dS$
in (\ref{eq:
weak N 2 B}).
Hence we see that (\ref{cond 5 B}) is valid
and
thus
\ba
(\Delta_{\tilde g}+k^2)\tilde u=\tilde f\quad\hbox{in
}N
\ea
in the sense of the Definition \ref{energysolution}.

So far,
we have proven that a (classical) solution to the Helmholtz
equation
on $M$ yields a finite energy solution
to the equation on $N$. Next,
we prove the converse, i.e.,  that
the Helmholtz equation on $N$
implies Helmholtz equation on $M$.

Assume that $\tilde u$ satisfies
Helmholtz equation (\ref{case 3 B})
on $(N,\tilde g)$ in the sense of
Definition \ref{energysolution}, with $\tilde
f\in L^2(N)$
supported
away from $\Sigma$. In particular, $\tilde u$
is  a measurable
function in $N$
satisfying (\ref{cond 1 A}),
(\ref{cond 2 A}),
and
(\ref{cond 3 A}).

Let  $u=\tilde u\circ F$ and $f=\tilde f\circ F$
on $M\setminus \gamma$.
Then we have
\beq
\label{analog27}
(\Delta_g+k^2)u_1=f_1=f|_{M_1 \setminus \gamma_1}
\quad \hbox{in }
M_1\setminus \gamma_1,\quad
u_1|_{\p
M_1}=h
\eeq
and
\beq \label{analog28}
(\Delta_g+k^2)u_2=f_2=f|_{M_2
\setminus \gamma_2} \quad \hbox{in }
M_2\setminus \gamma_2.
\eeq

By
conditions (\ref{cond 1 A}),
(\ref{cond 2 A}),
and (\ref{cond 3
A}),
we have that
\ba
& & |u_i|^2\in L^1(M_i\setminus
\gamma_i,|g|^{1/2}dx),\\
& & g_i^{jk}(\p_j u_i)(\p_k u_i)\in
L^1(M_i\setminus \gamma_i,|g|^{1/2}dx),
i=1,2.
\ea
Thus $u_i\in
H^1(M_i\setminus \gamma_i,dx)$.
As before, we see that then
\beq
\label{analog29}
(\Delta_g+k^2)u_1&=&f_1\quad \hbox{in }
M_1,\quad
u_1|_{\p M_1}=h,\\ \nonumber
(\Delta_g+k^2)u_2&=&f_2\quad
\hbox{in } M_2,
\eeq
where $f_i$ is extended to have the value $0$ at
$\gamma_i$ and $u_i$ are
smooth near $\gamma_i$.

Since $\tilde u$
satisfies (\ref{cond 5 B}),
we  see that for ${\tilde \psi} \in
C^\infty_0(N)$,
\ba
0&=&\int_{N} [-D_{\tilde g}\tilde u \p_j
\tilde
\psi+k^2 \tilde u\tilde \psi |\tilde g|^{1/2}-\tilde
f\tilde
\psi |\tilde g)|^{1/2}]dx
\\ \nonumber
&=&\lim_{\e\to 0}
\int_{N\setminus \Sigma(\e)} (-\tilde
g^{ij}\,\p_i \tilde u\, \p_j
\tilde \psi+(k^2 \tilde u+\tilde
f)\tilde \psi)\ |\tilde
g|^{1/2}dx\\
\nonumber
&=&\lim_{\e\to 0} \left(\int_{\p
\Sigma(\e)\cap N_2}+\int_{\p
\Sigma(\e)\cap N_1}
\right)
(-\tilde
g^{ij}\,\p_i \tilde u|_{\p \Sigma(\e)}\, \nu _j \tilde
\psi)
|\tilde
g_{bnd}|^{1/2}dS(x)
\nonumber\\
&=&
\lim_{\e\to 0} \int_{\p
\gamma_1(\e)}
(-g^{ij}\,\p_i u_1\, \nu _j  \psi)
|g_{bnd}|^{1/2}dS(x)
\\ \nonumber
& &+\lim_{\e\to 0} \int_{\p
\gamma_2(\e)}
(-g_s^{ij}\,\p_i u_2\, \nu _j  \psi)
|g_{bnd}|^{1/2}dS(x)\\
&=&  \nonumber
0,
\ea
where $\psi =\tilde \psi
\circ F$. Here as
in the proof of Proposition \ref{single coating
with Helmholtz},
we use  the fact that
$u_1$ is smooth function
implying that $\p_i u_1$ is bounded.

Thus, we have shown that
the
function $u$ is a classical solution on $M$
of
\beq\label{analog33}
(\Delta_{g}+k^2)u_1=f_1 \quad\hbox{in
}M_1,\quad
u_1|_{\p
M_1}=h
\eeq
and
\beq\label{analog34}
(\Delta_{g}+k^2)u_2=f_2
\quad\hbox{in }M_2.
\eeq
This proves the claim.
\hfill\proofbox

Next we prove a result that is not necessary for
the proof  but
gives, in the case of the double coating,
      an
alternative  treatment
of the distribution $D_{\tilde g}\tilde u$,
simpler than before.

\begin{lemma} In the double coating
construction, the term
\beq \label{cond 4}
|\tilde g|^{1/2}\tilde
g^{ij}\p_i\tilde u\in {\cal D}'(N,dx),
\eeq
appearing in
Definition
\ref{energysolution}  {as $D_{\tilde g}\tilde u$}
is
well-defined as a sum of products of
Sobolev distributions and
Lipschitz functions.
\end{lemma}

\noindent{\bf Proof.}
The problem
we need to consider is here is that
$ L^2(N,\, |\tilde g|^{1/2}dx)$
contains
functions that are not locally integrable with respect
to
measure $dx$ and thus we do not immediately see that
distribution
derivatives $\p_j\tilde u$ in $N$ are well
defined. We deal with this
by applying condition (\ref{cond 3 B}).
To do this, let $u=\tilde
u\circ F:M\setminus \gamma\to \R$.
Using (\ref{cond 1 B}), (\ref{cond
2 B}), (\ref{cond 3 B})
and changing variables in the integration,
one sees that
\ba
\int_{M\setminus\gamma} |g|^{1/2}g^{ij}(\p_i u)\p_j
u)\,dx<\infty.
\ea
As $g$ is bounded from above and below,
this
implies that $u\in H^1(M\setminus\gamma,dx)
\subset L^6
(M\setminus\gamma),dx)$. Furthermore, changing
variables again
implies that
\ba
\int_{N\setminus\Sigma} |\tilde g|^{1/2}|\tilde
u|^6\,dx
<\infty,
\ea
so that $\tilde u\in
      L^6(N,\, \det (\tilde
g)^{1/2}dx).$
Now in the  boundary normal coordinates $(\omega,\tau)$
near
$ \Sigma$,
$\tau(x)=\hbox{dist}_{\R^3}(x, \Sigma)$,
we
have
\ba
\tau^{-2}|\tilde g|^{1/2}\in [c_1,c_2],\quad
c_1,c_2>0,
\ea
and thus
\ba
\int_{N} |\tilde u|\,dx&=&
\int_{N}
|\tilde u|\tau(x)^{1/3} \tau(x)^{-1/3} \,dx\\
&\leq &
\|\tilde
u\,\tau^{1/3}
\|_{L^6(N,dx)}
\|\tau(x)^{-1/3}\|_{L^{6/5}(N,dx)}\\
&\leq &
\|\tilde
u \|_{L^6(N,\tau^2 dx)}
\|\tau(x)^{-2/5}\|_{L^1(N,dx)}\\
&\leq
&
C\|\tilde u \|_{L^6(N,|\tilde g|^{1/2} dx)}<\infty,
\ea
cf. the
discussion at the end of \S\S 3.2.3.
A similar computation shows that
$ \tilde u\in L^p(N,dx)$
for some $p>1$, and thus $\p_j \tilde u\in
W^{-1,p}(N,dx)$.
As is shown  at the end of \S 2 that
\beq\label{eq:
lip}
|\tilde g|^{1/2}\tilde g^{jk}\in C^{0,1}(N),
\eeq
multiplication
by $|\tilde g|^{1/2} \tilde g ^{jk}$ maps
$W^{1,p'}\longrightarrow
W^{1,p'}$ and thus, by duality,
\ba
|\tilde g|^{1/2}\tilde g^{jk}
\p_j \tilde u\in W^{-1,p}(N,dx),
\ea
i.e., the distribution
(\ref{cond 4})
is defined as a sum of products of Lipschitz functions
and
$W^{-1,p}$-distributions
\hfill\proofbox
\bigskip

\subsection{
Coating with a lining: a physical surface }

In the previous sections
we have considered the Helmholtz
equation in a domain $N\subset
\R^3$, equipped with a  metric $\tilde g$
that is singular at a
surface $\Sigma$.
Later, for Maxwell's equations, we will need to
consider
$\Sigma$ as a  ``physical" surface, i.e., an obstacle
on
which  we have to impose a boundary condition.
To motivate these
constructions, we consider next,
for the Helmholtz equation, what
happens when we have such a physical
surface at $\Sigma$. More
precisely,
we consider the Helmholtz equation
in the domain
$N\setminus \Sigma=N_1\cup N_2$
where, on the  both sides of the
boundary of $\Sigma$, that is,
on $\Sigma_+=\p N_1\setminus\p N$ and
on $\Sigma_-=\p N_2$,
we impose a degenerate boundary condition of Neumann type. In
physical terms,
this corresponds to having a material, sound
hard surface
located at $\Sigma$,
separating  space into two open
components, $N_1$ and $N_2$.
Although we will not need this, it can in fact be shown that $\tilde 
u$ is a solution
in the sense of Def. \ref{doubcoatingenergy N} iff it is in the sense of
Def. \ref{energysolution}.

\subsubsection{Weak solutions  for the double coating
with
Neumann boundary conditions}

In the following, we consider a
double coating $(M,N,F,\Sigma,g)$.
Suppose that $k\geq 0$ and $\tilde
f\in L^2 (N, |\tilde g|^{1/2}
dx)$.

\begin{definition}\label{doubcoatingenergy N}
We say that
$\tilde u$ is a \emph{finite energy} solution
of the boundary value
problem
with degenerate Neumann boundary conditions at $\Sigma$,
\beq\label{case
3 N}
& &(\Delta_{\tilde g}+k^2)\tilde u=\tilde f\quad\hbox{in
}N\setminus \Sigma,\\
& &
\tilde u|_{\p N}=\tilde h \label{case 3 Na}
\\  \label{case 3 Nb}
& &
|\tilde g|^{1/2} \p_\nu \tilde u|_{\Sigma_+}=0,
\quad 
|\tilde g|^{1/2}  \p_\nu
\tilde u|_{\Sigma_-}=0,
\eeq
if $\tilde u$
is a measurable function
in $N\setminus \Sigma $ such that
\beq\label{cond 1 N}
& &\tilde u\in
L^2(N\setminus \Sigma,\, |\tilde g|^{1/2}dx);\\
\label{cond 2 N}
&
&\p_j\tilde u\in H^1_{loc}(N\setminus \Sigma,dx);
\\
\label{cond 3
N}
& &\int_{N\setminus\Sigma} |\tilde g|^{1/2}\tilde
g^{ij}\p_i\tilde
u\p_j\tilde u\,dx<\infty;\\
\label{cond 4bb N}
&
&
(\Delta_{\tilde g}+k^2)\tilde u=\tilde f\hbox{ in some neighborhood
of }\p N,\\
& &\nonumber
\tilde u|_{\p N}=\tilde h;
\eeq
and finally,
\beq\label{cond 5 N}
\int_{N\setminus \Sigma} \left(-\tilde g^{ij}\,\p_i \tilde u\, \p_j
\tilde \psi+(k^2-\tilde f) \tilde u \tilde \psi\right)|\tilde g|^{1/2}
dx=0
\eeq
for all
\bfo
\tilde \psi=
\left\{
\begin{array}{cl}
\tilde \psi_1(x), & x \in N_1, \\
\tilde \psi_2(x), & x \in N_2,
\end{array}
\right.
\efo
with $\tilde \psi_1\in C^\infty(\overline N_1)$ vanishing
near the exterior boundary $\p N = \p N_1\setminus \Sigma$  and
        $\tilde \psi_2\in C^\infty(\overline N_2)$.
\end{definition}

{Invisibility  for the
double coating with a physical surface at $\Sigma$, with
respect to the Helmholtz equation  at all frequencies, is a consequence of the
following analogue of  Theorem
\ref{double coating with Helmholtz:a} :}

\begin{theorem}
\label{double coating with Helmholtz N}
Let $u=(u_1, u_2):M\setminus \gamma\to \R$
and  $\tilde u:N\setminus \Sigma\to \R$ be  measurable functions
such that $u=\tilde u\circ F$.
Let $f=(f_1, f_2):M\setminus \gamma\to \R$
and  $\tilde f:N\setminus \Sigma\to \R$ be $L^2$ functions
supported away from $\Sigma$ and $\gamma$
such that  $f=\tilde f\circ F$,
and $\tilde h:\p N\to\R,\, h:\p M_1\to\R$ be such that $h=\tilde
h\circ F_1$..

Then the following are equivalent:
\begin{enumerate}

\item The function
$\tilde u$, considered as a measurable function on $N\setminus \Sigma$,
is
a finite energy solution of (\ref{case 3 N})
with Neumann boundary conditions at $\Sigma$
and
inhomogeneity $\tilde f$ in the sense of Definition \ref{doubcoatingenergy N}.

\item The function $u$ satisfies
\beq\label{Disc N}
(\Delta_{g}+k^2)u_1=f_1\quad\hbox{ on }M_1,\quad
u|_{\p M_1}=h:= \tilde h \circ F
\eeq
and
\beq\label{Ball N}
& &(\Delta_{g}+k^2)u_2=f_2\quad\hbox{ on }M_2.
\eeq

\end{enumerate}

\end{theorem}

\noindent
{\bf Proof.}  The proof is identical to that  of
Theorem  \ref{double coating with Helmholtz:a}.\hfill\proofbox
\bigskip

{\bf Remark.} Let  $\tilde g$ be a singular metric on  $N$ corresponding to a
double coating. The implication of Theorems
\ref{double coating with Helmholtz:a}  and \ref{double coating with
Helmholtz N} is that the solutions $\tilde u$ in $N\setminus \Sigma$
coincide in the following
cases:
\begin{enumerate}

\item  We have the metric $\tilde g$ on $N$, singular at the virtual surface
$\Sigma$.

\item We have the metric  $\tilde g$ on $N\setminus \Sigma$
and a sound hard  physical surface at  $\Sigma$,
in the sense of Definition  \ref{doubcoatingenergy N}.

\end{enumerate}

Similar results can be proven  when
the metric $\tilde g$ in $N$ corresponds to  a single coating.

\section{Maxwell's equations}

\subsection{Geometry and definitions}

Let us start with a general Riemannian manifold $(M,g)$, possibly
with a  non-empty
boundary, and consider how to define  Maxwell's equations on $M$.
We follow the treatment in \cite{KLS}, using, however, slightly different
notation.

Using the  metric $g$,
we define a permittivity
and permeability by setting
\ba
& &\e^{jk}=\mu^{jk}=|g|^{1/2}g^{jk},\quad \hbox{on }M.
\ea

Although defined with respect to local coordinates, $\e$ and $\mu$
are in fact invariantly defined, and transform as a product of a $(+1)-$density
and a  contravariant symmetric two-tensor.
\medskip

{\bf Remark.}
In   $\R^3$
with the Euclidean metric $g_{jk}=\delta_{jk}$, we  have
$\e^{jk}=\mu^{jk}=\delta^{jk}$. If we would like to define a generalization of
isotropic media on a general Riemannian manifold, it would be as
\ba
& &\e^{jk}=\alpha(x)^{-1}|g|^{1/2}g^{jk},\\
& &\mu^{jk}=\alpha(x)|g|^{1/2}g^{jk},
\ea
on $M$, where $\alpha(x)$ is a positive scalar function. However, in
the following
we assume for simplicity that $\alpha=1$.

\medskip

In the following we consider the electric and magnetic
fields, $E$ and $H$, as
differential 1-forms, given in some local
coordinates by
\ba
E=E_jdx^ j,\quad H=H_jdx^ j,\quad
\ea
and $J$, the
internal current, as a 2-form.

Now consider the time harmonic  Maxwell's equations  on  $(M,g)$ at frequency
\nolinebreak$k$. They can  be written invariantly as
\beq\label{invariant maxwell}
dE=ik*_gH,\quad dH=-ik*_gE + J
\eeq
where $*_g:C^\infty(\Omega^jM)\longrightarrow C^\infty(\Omega^{3-j}M)$ denotes
the Hodge-operator on \linebreak$j$-forms, $0\le j\le 3$, given on 1-forms by
\beq\label{conc 1}
*_g(E_jdx^ j)&=&
\frac{1}{2}\vert g\vert^{1/2} g^{j l} E_j\,
s_{lpq}dx^p\wedge dx^q\\
&=&
\frac{1}{2}\e^{j l} E_j\,
s_{lpq}dx^p\wedge dx^q\nonumber
\eeq
where $s_{lpq}$ is the  Levi-Civita permutation symbol,
$s_{lpq}=1$ if $(l,p,q)$ even permutation of $(1,2,3)$,
$s_{lpq}=-1$ if $(l,p,q)$ odd permutation of $(1,2,3)$, and
zero otherwise. Thus
\ba
*_g(E_jdx^ j)=
(\e^{j3} E_j)\, dx^1\wedge dx^2-(\e^{j 2} E_j)\, dx^1\wedge dx^3+
(\e^{j 1} E_j)\,dx^2\wedge dx^3.
\ea

Next, we want to write these equations in arbitrary coordinates
so that they resemble the traditional Maxwell equations. The idea is
that
we want to have expressions which specialize, in the case of
the Euclidean metric on $\R^3$,  to expressions involving  $curl$ and
the matrices
$\e^{jk}$ and $\mu^{jk}$.
To write equations in such a  form,
let us introduce, for $H=H_jdx^j$,
the notation
\ba
(\hbox{curl}\,H)^l=s^{lpq}\frac {\p }{\p x^p}H_q\quad.
\ea
The exterior derivative
\ba
d(H_jdx^ j)&=& \frac {\p H_j}{\p x^k}\,dx^k\wedge dx^j
\ea
may then be written as
\beq\label{conc 2}
dH&=&\frac 12 (\hbox{curl}\,H)^l\,s_{lpq}dx^p\wedge dx^q.
\eeq

Combining (\ref{conc 1}) and (\ref{conc 2}) we see
that Maxwell equations (\ref{invariant maxwell}) can be written as
\ba
& &(\hbox{curl}\,E)^l=ik\,\mu^{j l} H_j,\\
& &(\hbox{curl}\,H)^l=-ik\,\e^{j l} E_j+J^l.
\ea

Below, we denote also
\ba
(\nabla \times E)^j=(\hbox{curl}\,E)^j,
\ea
and usually denote the standard volume element of $\R^3$ by $dV_0(x)$.

There are many boundary conditions that makes the boundary
value problem for Maxwell's equations on a domain well posed. For example:

\begin{itemize}

\item Electric boundary condition:
\ba
\nu\times E|_{\p M}=0,
\ea
where $\nu$ is the Euclidean normal vector of $\p M$.
Physically this corresponds to lining the boundary with a  perfectly
conducting material.

\item Magnetic boundary condition:
\ba
\nu\times H|_{\p M}=0,
\ea
where $\nu$ is the Euclidean normal vector of $\p M$.
In other
words, the tangential components of the magnetic field vanish.

\item
Soft and hard surface (SHS) boundary condition
\cite{HLS,Ki,Ki2,Li}:
\ba
\zeta\,\cdotp E|_{\p M}=0\quad\hbox{and}\quad
\zeta\,\cdotp H|_{\p M}=0
\ea
where $\zeta=\zeta(x)$ is a tangential
vector field on $\p M$,
that is, $\zeta\times \nu=0$.
In other words,
the part of the  tangential component of the electric
field $E$ that
is parallel to $\zeta$ vanishes, and
the same is true for the
magnetic field $H$.
This can be physically realized by having a
surface with
thin parallel gratings
\cite{HLS,Ki,Ki2,Li}.
\end{itemize}

\subsection{Definition of
solutions of Maxwell equations}

Assume that $k\in \R\setminus
\{0\}$. We will define
finite energy solutions for
Maxwell's
equations in the same way for both the
single and double
coatings.

Let $(M,N,F,\gamma,\Sigma,g)$ be either a single or
double
coating construction, as in \S\ref{gbc},
denoting as usual $\tilde g=F_* g$ on  $ N\setminus \Sigma$.
On $M$ and $ N\setminus \Sigma$, we then define permittivity
and permeability tensors by setting
\ba
& &\e^{jk}=\mu^{jk}=|g|^{1/2}g^{jk},\quad \hbox{on }M,\\
& &\tilde \e^{jk}=\tilde \mu^{jk}=|\tilde g|^{1/2}\tilde
g^{jk},\quad \hbox{on }
            N\setminus \Sigma.
\ea
Let $J$ be a smooth internal current 2-form on $M$ that is supported away
from $\p M$.

\subsection{Finite energy solutions for single \\  and double coatings}

The definition  of finite energy solution is the same for both coatings.
On $M$, the parameters $\e$ and $\mu$ are bounded from
below and above, so  Maxwell's equations,
\beq\label{eq: physical Max a}
& &\nabla\times  E = ik  \mu(x)   H,\quad \nabla\times
            H =-ik  \e(x) E+J\quad
\hbox{ in }M,\\
& &R(\nu, E, H)|_{\p M}=b \nonumber
\eeq
are defined in the  sense of distributions in the usual way.
Here, $\nu$ denotes the Euclidian unit normal vector
of $\p M$ and $R(\cdot,\cdot,\cdot)$ is a boundary value operator
corresponding to the boundary conditions of interest, e.g.,
$R(\nu, E, H)=\nu\times E$ for the electric boundary condition.

If $J$ is smooth,  Maxwell's equations  imply that $E,H\in C^\infty
(M)$.

Next, we consider Maxwell's equations on $N$.
Let $\tilde J$ be a smooth 2-form on $N$ that is supported away
form $\p N\cup \Sigma$.

\begin{definition}\label{Maxwell-def}
Let $(M,N,F,\gamma,\Sigma,g)$ be either a single or double
coating.
We say that  $(\tilde E, \tilde H)$ is a finite energy solution
to Maxwell's equations on $N$,
\begin{equation}\label{eqn-4.1-main}
\nabla\times \tilde E = ik \tilde \mu(x)  \tilde H,\quad \nabla\times
            \tilde H =-ik  \tilde \e(x) \tilde E+\tilde J\quad
\hbox{ on }N,
\end{equation}
if $ \tilde E$, $\tilde H$, $\tilde D:=\tilde \e\, \tilde E$ and
   $\tilde B:=\tilde \mu\, \tilde H$
are forms in $N$ with $L^1(N,dx)$-coefficients satisfying
\beq\label{eq: Max norm}
\|\tilde E\|_{L^2(N,|\tilde g|^{1/2}dV_0(x))}^2= \int_{N}
\tilde \e^{ij}\, \tilde E_j\, \overline{  \tilde E_k} \,dV_0(x)<\infty
,\\
\|\tilde H\|_{L^2(N,|\tilde g|^{1/2}dV_0(x))}^2=\int_{N}
\tilde \mu^{ij}\, \tilde H_j\, \overline{ \tilde H_k} \,dV_0(x)<\infty;
\eeq
$(\tilde E,\tilde H)$ is a classical solution of Maxwell's equations
on a neighborhood
$U\subset \overline N$ of $\p N$:
\ba
& &\nabla\times  \tilde E = ik  \tilde \mu(x)   \tilde H,\quad \nabla\times
\tilde  H =-ik  \e(x)\tilde  E+ \tilde J\quad
\hbox{in }U,\\
& &R(\nu, \tilde E,\tilde H)|_{\p N}=\tilde b; \nonumber
\ea
and finally,
\ba
& &\int_N ((\nabla\times \tilde h)\,\cdotp \tilde E-
             ik \tilde h \,\cdotp \tilde \mu(x)\tilde H) \,dV_0(x)=0,\\
& &
\int_N (
(\nabla\times \tilde e)\,\cdotp \tilde H+
             \tilde e \,\cdotp ( ik\tilde \e(x)\tilde E-\tilde J)) \,dV_0(x)=0
\ea
for all $\tilde e,\tilde h\in C^\infty_0(\Omega^1N)$.
\end{definition}
Here, $ C^\infty_0(\Omega^1N)$ denotes smooth 1-forms on $N$
whose supports do not intersect $\p N$,
and the inner product ``$\cdotp$'' denotes the Euclidean inner product.

{\bf Remark.} The fact that $\tilde E,\tilde H$ are solutions of 
(\ref{eqn-4.1-main}) in the
sense of Def. \ref{Maxwell-def} implies that they are distributional 
solutionsin the usual
sense. Thus they also satisfy  the divergence equations,
\beq\label{eq: div 1B}
\nabla\cdotp \tilde \e \tilde E =\frac 1{ik}\nabla\cdotp \tilde J
,\quad \nabla\cdotp\tilde  \mu\tilde  H =0,
\eeq
in the sense of distributions.

\section{Full wave invisibility for the double coating}

In this section,  $(M,N,F,\gamma,\Sigma,g)$ denotes
a double coating construction. Invisibility for active devices
enclosed in the double coating,
with respect to Maxwell's equations at all frequencies, is a consequence of:

\begin{theorem}\label{double coating with Maxwell}
Let $E$ and $H$ be 1-forms
with measurable coefficients
on $M\setminus \gamma$ and $\tilde E$ and $\tilde H$ be
1-forms
with measurable coefficients on $N\setminus \Sigma$
such that $E=F^*\tilde E$, $H=F^*\tilde H$.
Let $J$ and $\tilde J$ be  2-forms
with smooth coefficients on $M\setminus \gamma$ and
$N\setminus \Sigma$ that are supported away
from $\gamma$ and $\Sigma$.

Then the following are equivalent:
\begin{enumerate}

\item
The
1-forms $ \tilde E$ and $ \tilde H$ on $N$ form a finite energy solution of
Maxwell's equations
\beq\label{eq: physical Max}
& &\nabla\times \tilde E = ik \tilde \mu(x)  \tilde H,\quad \nabla\times
           \tilde H =-ik  \tilde \e(x) \tilde E+\tilde J\quad
\hbox{ on }N,
\\
& &\nonumber R(\nu, \tilde E,\tilde H)|_{\p N}=b.
\eeq

\item
The 1-forms $E$ and $H$ on $M$ satisfy Maxwell's equations
\ba
& &\nabla\times   E = ik   \mu(x)    H,\quad \nabla\times
             H =-ik    \e(x)   E+  J\quad
\hbox{ on }M_1,\\
& &\nonumber R(\nu, E, H)|_{\p N}=b
\ea
and
\ba
& &\nabla\times   E = ik   \mu(x)    H,\quad \nabla\times
             H =-ik    \e(x)   E+  J\quad
\hbox{ on }M_2.
\ea

\end{enumerate}
\end{theorem}

\noindent
{\bf Proof.}
First we prove that Maxwell's equations
on $M$ imply Maxwell equations on $N$

Assume now that the 1-forms
$E$ and $H$ are classical solutions of
Maxwell's equations
on $M=M_1\cup M_2$,
\beq\nonumber
& &\nabla\times    E = ik    \mu(x)     H,\quad \nabla\times
              H =-ik  \e(x)   E + J\quad
\hbox{ on }M=M_1\cup M_2,\\
& &R(\nu, E, H)|_{\p N}=b. \label{eq: physical Max AAA}
\eeq
Since $J$ vanishes near $\gamma$, ellipticity implies  that
$E$ and $H$ are smooth near $\gamma$.

Define on $N\setminus \Sigma$ the forms
$\tilde E=(F^{-1})^* E,$
$\tilde H=(F^{-1})^* H,$
and $\tilde J=(F^{-1})^* J.$

Then $\tilde E$ satisfies the Maxwell's equations on
$N\setminus
\Sigma$,
\beq\label{eq: physical Max2}
\nabla\times \tilde E = ik
\tilde \mu(x)  \tilde H,\quad \nabla\times
           \tilde H =-ik
\tilde \e(x) \tilde E+\tilde J\quad
\hbox{ on }N\setminus
\Sigma,
\eeq

Again, let $\Sigma(t)$ be the $t$-neighborhood
of
$\Sigma$ with respect to the metric
$\tilde g$ and
$\gamma(t)$  the
$t$-neighborhood
of $\gamma$
with respect to $g$.
Let $I_t:\p \gamma
(t)\to M$ be the identity embedding.
Denote by $\nu$ be the
unit normal vector of  $\p\Sigma(t)$
and $\p\gamma(t)$ in
Euclidean metric.

Now, writing $E=E_j(x)dx^j$ on $M$, we see
using
the transformation rule for differential
1-forms that
the form
$\tilde E=(F^{-1})^*E$ is in local
coordinates is
\ba
\tilde E=\tilde
E_j(\tilde x)d\tilde x^j=(DF^{-1})_j^k(\tilde
x)\,E_k(F^{-1}(\tilde
x))
d\tilde x^j
,\quad
\tilde x\in N\setminus \Sigma,
\ea
and, using
$F_t= F\circ
I_t:\p \gamma(t)\to \p \Sigma(t)$, we
have
\beq\label{tilde I formula}
\tilde
I^*(\tilde
E_j(x)dx^j)=(DF_t^{-1})_j^k(\tilde x)\,E_k(F^{-1}(\tilde
x))\,
d\tilde x^j,\quad \tilde x=F(x)
\eeq
Let us now do
computations
in the Euclidean coordinates.
In the Euclidean metric $g_e$,
the
matrix  $DF_t^{-1}$
satisfies
\beq\label{transformation
1}
\|DF_t^{-1}\|_{(T\p\Sigma(t),g_e)\to
(T\p\gamma(t),g_e)}\leq
Ct,
\eeq
and
since $E$ is smooth  near $\gamma$ we see
\ba
|\nu\times \tilde
E(y)|_{\R^3}\leq Ct,\quad
y\in\p\Sigma(t).
\ea
Thus using (\ref{eq:
physical Max2})
we see that
for $\tilde h\in
C^{\infty}_0(\Omega^1N)$
\beq\label{eq: weak N 2
Max}
&
&\int_N
((\nabla\times \tilde h)\,\cdotp \tilde E-
            ik
\tilde h \,\cdotp \tilde \mu(x)\tilde H) \,dV_0(x)\\
& &=\lim_{t\to
0} \int_{N\setminus \Sigma(t)}
((\nabla\times \tilde h)\,\cdotp
\tilde E-
            ik \tilde h \,\cdotp \tilde \mu(x)\tilde H)
\,dV_0(x)
\nonumber
\\
& &=-\lim_{t\to 0} \int_{\p \Sigma(t)}
(\nu
\times \tilde E)\,\cdotp \tilde h\,  dS(x)=0. \nonumber
\eeq

Thus,
we have shown that
\beq\label{eq: physical Max3A}
\nabla\times \tilde
E = ik \tilde \mu(x)  \tilde H\quad \hbox{ in }N
\eeq
in the sense of
Definition \ref{Maxwell-def}. Similarly, we see that
\beq\label{eq:
physical Max3B}
           \nabla\times  \tilde H =-ik  \tilde \e(x)
\tilde E+\tilde J
\quad \hbox{ in }N
\eeq
in the same finite energy
sense.

Next we show that Maxwell's equations on $N$ implies
Maxwell's
equations on $M$.
Let $U\subset M$ be a
bounded
neighborhood
of $\gamma$ and $W=F(U\setminus \gamma)\cup
\Sigma$ be
a
neighborhood of $\Sigma$ such that $\hbox{supp}\,(\tilde
J)\cap
W=
\emptyset$.

Assume that $\tilde E$ and $\tilde H$ form a
finite energy solution of
Maxwell's equations (\ref{eq: physical
Max})
on $(N,g)$ in finite energy sense
with a source $\tilde J$
supported
away from $\Sigma$, implying in
particular that
\ba
\tilde
\e^{jk}\tilde E_j\overline {\tilde
E_k} \in L^1(W,\,dx),\quad
\tilde \mu^{jk}\tilde H_j\overline
{\tilde H_k} \in
L^1(W,\,dx).
\ea
Define  $E=F^*\tilde E$,
$H=F^*\tilde H$ and
$J=F^*\tilde J$ on $M\setminus \gamma$.  We
have
\ba
\nabla\times E =
ik \mu(x) H,\quad \nabla\times H =-ik \e(x)
E+J
\quad \hbox{in }
M\setminus \gamma
\ea
and
\ba
      \e^{jk}  E_j
\overline {E_k} \in
L^1(U\setminus \gamma,\,dV_0(x)),\quad
      \mu^{jk}
H_j  \overline
{H_k} \in L^1(U\setminus \gamma,\,dV_0(x)).
\ea
As $\e$ and $\mu$ on
$M$ are bounded from above and below, these
imply that
\ba
& &
\nabla\times E \in
L^2(U\setminus
\gamma,\,dV_0(x)),\quad
\nabla\times H \in
L^2(U\setminus
\gamma,\,dV_0(x)),\\
& &\nabla\cdotp (\e
E)=0,\quad
\nabla\cdotp (\mu H)=0\quad \hbox{in } U\setminus
\gamma.
\ea
Let
$E^e,H^e\in L^2(U,\,dV_0(x))$ be measurable extensions
of $E$ and $H$
to $\gamma$.
Then
\ba
           & &\nabla\times E^e
-ik \mu(x) H^e=0\quad \hbox{in }
U\setminus \gamma,\\
           &
&\nabla\times E^e -ik \mu(x) H^e\in H^{-1}(U,\,dV_0(x)),\\
  & &\nabla\times H^e +ik \e(x) E^e=0\quad \hbox{in }
U\setminus
\gamma,\\
           & &\nabla\times H^e +ik \e(x) E^e\in
H^{-1}(U,\,dV_0(x)).
\ea
Since  $\gamma$ is a subset with (Hausdorff)
dimension 1
of the 3-dimensional domain $U$, it has zero capacitance.
Thus, the Lipschitz functions on $U$ that vanish on $\gamma$
are dense in $H^1(U)$, see \cite[Thm 4.8 and remark 4.2(4)]{KKM},
or \cite[Thm. 3.28] {AF}.
Thus there  are no non-zero  distributions in $H^{-1}(U)$ supported on
$\gamma$. Hence
        we see that
\ba
\nabla\times E^e -ik \mu(x) H^e=0,\quad
\nabla\times H^e +ik \e(x) E^e=0\quad \hbox{in } U.
\ea
This also implies that
\ba
\nabla\cdotp (\e E^e)=0,\quad \nabla\cdotp (\mu H^e)=0\quad \hbox{in } U.
\ea
These imply that $E^e$ and $H^e$ are in $C^\infty$ smooth in $U$.

Summarizing, $E$ and $H$ have unique continuous extensions
to $\gamma$, and the extensions are classical
solutions to Maxwell's equations.

\section{ Cauchy data for the single coating must vanish}

In this section  $(M,N,F,\gamma,\Sigma,g)$ denotes a single coating
construction.
The following gives
the counterpart for Maxwell's equations of the hidden Neumann boundary
condition  on $\p M_2$ that appeared for the Helmholtz equation.

\begin{theorem}\label{single coating with Maxwell A}
Let $E$ and $H$ be 1-forms
with measurable coefficients
on $M\setminus \gamma$ and $\tilde E$ and $\tilde H$ be
           1-forms
with measurable coefficients on $N\setminus \Sigma$
such that $E=F^*\tilde E$, $H=F^*\tilde H$.
Let $J$ and $\tilde J$ be  2-forms
with smooth coefficients on $M\setminus \gamma$ and
           $N\setminus \Sigma$, that are supported away
from $\gamma$ and $\Sigma$.

Then the following are equivalent:
\begin{enumerate}

\item
The 1-forms $ \tilde E$ and $ \tilde H$ on $N$ satisfy Maxwell's equations
\beq\label{eq: physical Max M}
& &\nabla\times \tilde E = ik \tilde \mu(x)  \tilde H,\quad \nabla\times
           \tilde H =-ik  \tilde \e(x) \tilde E+\tilde J\quad
\hbox{ on }N,
\\
& &\nonumber \nu\times \tilde E|_{\p N}=f
\eeq
in the sense of Definition \ref{Maxwell-def}.

\item
The forms $E$ and $H$ satisfy Maxwell's equations on $M$,
\beq\label{eq: physical Max single M1 new}
& &\nabla\times   E = ik   \mu(x)    H,\quad \nabla\times
             H =-ik    \e(x)   E+  J\quad
\hbox{ on }M_1,
\\
& &\nonumber \nu\times E|_{\p M_1}=f
\eeq
and
\beq\label{eq: physical Max single M2 new}
& &\nabla\times   E = ik   \mu(x)    H,\quad \nabla\times
             H =-ik    \e(x)   E+  J\quad
\hbox{ on }M_2
\eeq
with Cauchy data
\beq\label{eq: physical Max single M3 new}
& &\nu\times   E|_{\p M_2}=b^e,\quad
\nu\times   H|_{\p M_2}=b^h
\eeq
that satisfies $b^e=b^h=0$.

\end{enumerate}

Moreover, if $E$ and $H$ solve (\ref{eq: physical Max single M1 new}),
(\ref{eq: physical Max single M2 new}), and (\ref{eq: physical Max
single M3 new})
with non-zero $b^e$ or $b^h$, then the fields
$\tilde E$ and $\tilde H$ are not solutions of Maxwell equations
on $N$ in the sense of Definition \ref{Maxwell-def}.
\end{theorem}

\noindent {\bf Proof.}
Assume first that the 1-forms
$E$ and $H$ are classical solutions of
Maxwell's equations in $M$.
Moreover, assume that both $E$ and $H$ satisfy homogeneous
boundary condition
\beq\label{electric and magnetic condition}
\nu\times E|_{\p M_2}=0,\quad \nu\times H|_{\p M_2}=0,
\eeq
that is, for the field in $M_2$ the Cauchy data on $\p M_2$
vanishes. (Here,  $\nu$ again denotes the  Euclidean
unit normal of these surfaces.)

Again, define on $N\setminus \Sigma$ forms
$\tilde E=(F^{-1})^* E,$
$\tilde H=(F^{-1})^* H,$ and
$\tilde J=(F^{-1})^* J$.
Then $\tilde E$ satisfies  Maxwell's equations
on
$N\setminus \Sigma$,
\beq\label{eq: physical Max2 C}
\nabla\times
\tilde E = ik \tilde \mu(x)  \tilde H,\quad \nabla\times
           \tilde H
=-ik  \tilde \e(x) \tilde E+\tilde J\quad
\hbox{ in }N\setminus
\Sigma,
\eeq

Again, let $\Sigma(t)$ be the $t$-neighborhood
of
$\Sigma$ in the  $\tilde g$-metric and
$\gamma(t)$ be the
$t$-neighborhood
of $\gamma$ in the $g$-metric.

Arguing as in (\ref{transformation 1}) and below,
we see that
\beq\label{E vanishes
on boundary}
|\nu\times \tilde E(y)|_{\R^3}\leq Ct,\quad
y\in\p\Sigma(t)\cap N_2.
\eeq

Recall that $\Sigma_1(\e)=N_1\cap\Sigma(\e)$.
Then, using (\ref{eq: physical
Max2})
we see that for $\tilde h\in
C^\infty_0(\Omega^1N)$,
\beq\label{eq: weak N 2 Max b}
&
&\int_N
((\nabla\times \tilde h)\,\cdotp \tilde E-
            ik \tilde h
\,\cdotp \tilde \mu(x)\tilde H) \,dV_0(x)\\
&=&\lim_{t\to 0}
\int_{(N\setminus \Sigma_1(t)}
((\nabla\times \tilde h)\,\cdotp \tilde
E-
           ik \tilde h \,\cdotp \tilde \mu(x)\tilde
H)
\,dV_0(x)
\nonumber
\\
&=&-\lim_{t\to 0} \int_{\p
\Sigma_1(t)}
(\nu \times \tilde E)\,\cdotp \tilde h\,
dS(x)-
\int_{\p M_2}
(\nu \times \tilde
E)\,\cdotp \tilde h\,
dS(x)=0
\nonumber
\eeq
where
we used
(\ref{E vanishes on boundary})
and (\ref{electric and magnetic condition}).

Thus, we have shown
that
\beq\label{eq: physical Max3A 2}
\nabla\times \tilde E = ik
\tilde \mu(x)  \tilde H\quad \hbox{ on }N
\eeq
in the  sense of
Definition \ref{Maxwell-def}. Similarly, we see that
\beq\label{eq:
physical Max3B 2}
           \nabla\times  \tilde H =-ik  \tilde \e(x)
\tilde E+\tilde J
\quad \hbox{ on }N,
\eeq
also in the sense of
Definition \ref{Maxwell-def}.

Next we show that Maxwell's
equations on $N$ imply
Maxwell's equations on $M$.

Assume that
$\tilde E$  and $\tilde H$ form a finite energy solution of
Maxwell's
equations (\ref{eq: physical Max M}) on $(N,g)$.
Again, define on
$M\setminus \gamma$ forms
$E=F^*\tilde  E,$ $H=F^*\tilde  H$,
and
$J=F^*\tilde  J$.

As before, we see that $E$ and $H$ satisfy
Maxwell's equations
on $M_1\setminus
\gamma_1$ and the $E$ and $H$
are in $L^2(M_1,dV_0(x))$. Using
the removable of
singularity
arguments as in the case of double
coating,
we see that
$E$ and $H$ have extensions $E^e$
and $H^e$ in
$M_1$ that are classical solutions of
\beq
\label{ex1} & &\nabla\times E^e -ik \mu(x)
H^e=0\quad \hbox{ on } M_1,\\
\label{ex2} & &\nabla\times H^e +ik
\e(x) E^e=J\quad \hbox{ on } M_1.
\eeq
Note that (\ref{ex1}) implies that,
for the original field $\tilde E$,
\beq
\label{ex3}
\lim_{t\to 0}
\int_{\p \Sigma(t)\cap N_1}
(\nu \times \tilde E)\,\cdotp \tilde h\,
dS(x)=\lim_{t\to 0}
\int_{\p \gamma(t)\cap M_1}
(\nu \times E)\,\cdotp h\,
dS(x)=0
\eeq
where $h=F^*\tilde h$.

Moreover,   Maxwell's equations hold in the interior of $M_2$:
\ba
          \nabla\times E -ik \mu(x) H=0,\quad \nabla\times H +ik \e(x)
E=J\quad \hbox{ on }
M_2.
\ea

        Let us start to analyze what the validity of the equation
$\nabla\times \tilde
E -ik \tilde \mu(x) \tilde H=0$ on $N$
in the sense of Definition \ref{Maxwell-def} implies about the
boundary values on $\p M_2$.
Using (\ref{ex3}),
we see that for
$\tilde h\in C^\infty_0(\Omega^1N)$
\beq\label{eq: weak N 2 Max
C}
0&=&\int_N
((\nabla\times \tilde h)\,\cdotp \tilde E-
            ik \tilde
h \,\cdotp \tilde \mu(x)\tilde H) \,dV_0(x)\\
&=&\lim_{t\to 0}
\int_{(N\setminus \Sigma_1(t)}
((\nabla\times \tilde h)\,\cdotp \tilde E-
             ik \tilde h \,\cdotp \tilde \mu(x)\tilde H) \,dV_0(x)
\nonumber
\\
&=&-\left[\lim_{t\to 0} \int_{\p \Sigma_1(t)}
(\nu \times \tilde E)\,\cdotp \tilde h\,  dS(x)+\int_{\p N_2}
(\nu \times \tilde E)\,\cdotp \tilde h\,  dS(x)\right]
\nonumber
\\
&=&0- \int_{\p N_2}
(\nu \times E)\,\cdotp \tilde h\,  dS(x).
\eeq
This shows $
\nu\times E|_{\p M_2}=0.$
Similarly, the equation
$\nabla\times \tilde H +ik \tilde \e(x) \tilde E=\tilde J$ holding on $N$
in the finite energy sense  implies that
$
\nu\times H|_{\p M_2}=0.
$\hfill\proofbox
\bigskip

Assume that $E$ and $H$ satisfy the  time-harmonic
Maxwell's equations on $M_2\subset \R^3$ such that
the Cauchy data $(\nu\times E|_{\p M_2},\nu\times H|_{\p M_2})$
vanishes. By continuing $E$ and $H$ by zero to $\R^3\setminus M_2$
we obtain solutions of Maxwell's equation in $\R^3$.
Thus $J$ must  be a current for which there exist solutions of 
Maxwell's equations
in $\R^3$ both satisfying the Sommerfeld radiation condition  and
vanishing outside $N_2$. Such currents are nowhere dense in
$L^2(N_2)$, as
then the fields $E$ and $H$ corresponding to $J$ satisfy the
Sommerfeld radiation condition and,
using Stokes' theorem, we see that
the source $J$
is  orthogonal to all (vector-valued) Green's functions 
$G_e(\cdotp,y,k;a)$ with
$y\in \R^3\setminus
\overline M_2$ and $a\in \R^3$. Here, the Green's function
$(G_e(\cdotp,y,k;a),G_h(\cdotp,y,k;a))$ satisfies Maxwell's
equations in $\R^3$ with current $a\delta_y$
and the Sommerfeld radiation condition.

We thus conclude that finite energy  solutions to
Maxwell's equations on $N$ with
the single coating exist only if
the Cauchy data $(\nu\times E|_{\p M_2},\nu\times H|_{\p M_2})$
vanishes on the inner surface of the cloaked region.
Thus, finite energy solutions
do not exist for generic sources, i.e., internal currents $J$, in the cloaked
region.
\medskip

\section{Cloaking an infinite cylindrical domain}\label{cylinder}

We now consider an infinite cylindrical domain, $N=B_2(0,2)\times \R$ for
simplicity,
with the double coating. Here, $B_2(0,r)\subset \R^2$
is Euclidian disc with center $0$ and radius $r$.
Numerics for cloaking an infinite cylinder have been presented in
\cite{CPSSP}. {\newtext This may also provide a picture of the
cloaking that was physically
implemented with a ``sliced cylinder" geometry in \cite{SMJCPSS},
although precise modelling has
not been carried out. With this limitation in mind, the physical
interpretation of Theorems
\ref{single coating with Maxwell} and ref\{single coating with Maxwell
obstacle} below is that the
cloaking would be more effective with the insertion of a liner to
implement the SHS boundary
conditions which are necessary for the existence of finite energy solutions.}

Here, we modify the treatment from \S\ref{gbc} to the noncompact
setting,  blowing
up  a  line and trying to obtain
an infinitely long, invisible cable.

Let
\ba
& &M_1=B_2(0,2)\times \R,\quad
\gamma_1=\{(0,0)\}\times \R\subset M_1,\\
& &M_2=S^2\times \R,\quad
\gamma_2=\{NP\}\times \R\subset M_2
\ea
Let $M =M_1\cup M_2$, $\gamma=\gamma_1\cup \gamma_2$,
\ba
& &N_1=B_2(0,2)\times \R \setminus(\overline B_2(0,1)\times \R),\\
& &N_2=B_2(0,1)\times \R,\\
& &\Sigma=\p B_2(0,1)\times \R,
\ea
and $N=B_2(0,2)\times \R=N_1\cup N_2\cup \Sigma$.
Let
\ba
F=(F_1,F_2)&:&M\setminus \gamma\to N\setminus \Sigma
\ea
be such that
\ba
F_1&:&M_1\setminus \gamma_1\to
N_1,\\
F_2&:&M_2\setminus \gamma_2\to N_2.
\ea
are
diffeomorphisms.
Let $X:B_2(0,2)\times \R\setminus \{(0,0)\}\times
\R)\to (r,\theta,z)$
be the standard cylindrical coordinates on
$M_1$.
We assume that $F$ is stretching only in radial direction,
that
is,
\beq\label{rotations}
X(F(X^{-1}(r,\theta,z)))=(F_1(r),\theta,z).
\eeq
Similarly,
on $M_2$ we have
variables $(r,\theta,z)$, where $r=\dist(x,SP)$
and
we assume that $F$ has a form analogous to (\ref{rotations})
in
$M_2$.
For simplicity, let $g_1$ be the Euclidean metric on $M_1$
and $g_2$
the product
of standard metric
on $S^2$ and standard metric
of $\R$ on $M_2$.
Let
$\tilde g=F_* g$ on $N\setminus \Sigma$, so
that $(M,N,F,\gamma,\Sigma,g)$ is
a double coating construction in
this context.

On $M$ and $ N\setminus
\Sigma$ we define
permittivity
and permeability by
setting
\ba
&
&\e^{jk}=\mu^{jk}=|g|^{1/2}g^{jk},\quad \hbox{on
}M_1\cup
M_2,\\
& &\tilde \e^{jk}=\tilde \mu^{jk}=|\tilde
g|^{1/2}\tilde
g^{jk},\quad \hbox{on }
          N\setminus
\Sigma.
\ea
By finite energy
solutions of Maxwell's equations on $N$
we will mean one-forms
$\tilde E$ and $\tilde H$ satisfying
the conditions of Definition
4.1, where we emphasize the assumption
that $\tilde D:=\tilde \e \tilde E$ and
$\tilde B:=\tilde \mu \tilde H$
are in $L^1(N\setminus \Sigma,dx)$,
making the integrals
at the end of Definition 4.1 well defined.

To
formulate the results, we need to definethe restrictions
of fields on
the lines $\gamma_1\subset M_1$ and $\gamma_2\subset M_2$.
First,
assume  that the 1-forms
$E$ and $H$ on $M$ are classical solutions
to
Maxwell's equations
on $M$,
\beq\label{eq: physical Max BBB}
&
&\nabla\times    E = ik    \mu(x)   H,\quad
\hbox{in }M=M_1\cup
M_2,\\
\nonumber
& &\nabla\times H =-ik  \e(x)   E + J,\quad
\hbox{in
}M=M_1\cup M_2,
\\
& &\nu\times
E|_{\p M_1}=f, \nonumber
\eeq
where
$J$ is supported away from
$\gamma=\gamma_1\cup \gamma_2$.
Note that
then $E$ and $H$ are $C^\infty$ near
$\gamma$, and thus we can define
the restrictions
of the vertical
components of
the fields on $\gamma_1\subset
M_1$,
\beq\label{special electric and magnetic condition
1}
\zeta\,\cdotp E|_{\gamma_1}=b^e_1,\quad \zeta\,\cdotp
H|_{\gamma_1}=b^h_1,
\eeq
where $\zeta:=(0,0,1)=\frac \p {\p z},\, z:=x^3$.

Similarly, we can define
$b^e_2$ and
$b^h_2$ to be the restrictions on $\gamma_2\subset
M_2$,
\beq\label{special electric and magnetic condition
2}
\zeta\,\cdotp E|_{\gamma_2}=b^e_2,\quad \zeta\,\cdotp
H|_{\gamma_2}=b^h_2.
\eeq
Note that $b^e_j=b^e_j(z)$ and
$b^h_j=b^h_j(z)$, $j=1,2$, depend
only on $z$.

\begin{theorem}\label{single coating with Maxwell}
Let $E$ and $H$ be 1-forms
with measurable coefficients
on $M\setminus \gamma$ and $\tilde E$ and $\tilde H$ be
1-forms
with measurable coefficients on $N\setminus \Sigma$
such that $E=F^*\tilde E$, $H=F^*\tilde H$.
Let $J$ and $\tilde J$ be  2-forms
with smooth coefficients on $M\setminus \gamma$ and
$N\setminus \Sigma$, that are supported away
from $\gamma$ and $\Sigma$, respectively

Then the following are equivalent:
\begin{enumerate}

\item
On $N$, the
1-forms $ \tilde E$ and $ \tilde H$ satisfy Maxwell's equations
\beq\label{eq: physical Max tube}
& &\nabla\times \tilde E = ik \tilde \mu(x)  \tilde H,\quad \nabla\times
           \tilde H =-ik  \tilde \e(x) \tilde E+\tilde J\quad
\hbox{ in }N,
\\
& &\nonumber \nu\times \tilde E|_{\p N}=f
\eeq
and  $\tilde E$ and $\tilde H$ are finite energy
solutions.

\item
On $M$, the forms $E$ and $H$ are classical solutions to Maxwell's equations
(\ref{eq: physical Max BBB})
on $M$, with data
\beq\label{eq: M restrictions}
b^e_1=\zeta\cdotp E|_{\gamma_1},\quad
b^e_2=\zeta\cdotp E|_{\gamma_2},\quad
b^h_1=\zeta\cdotp H|_{\gamma_1},\quad
b^h_2=\zeta\cdotp H|_{\gamma_2},
\eeq
that satisfy
\beq\label{eq: M restrictions same}
b^e_1(z)=b^e_2(z)\quad\hbox{and}\quad b^h_1(z)=b^h_2(z),\quad  z\in \R.
\eeq

\end{enumerate}
Moreover, if $E$ and $H$ solve (\ref{eq: physical Max BBB})
with restrictions (\ref{eq: M restrictions}) that do not
satisfy (\ref{eq: M restrictions same}), then then the fields
$\tilde E$ and $\tilde H$ are not finite energy solutions of Maxwell equations
on $N$.
\end{theorem}

\noindent
{\bf Proof.}
First we show that the
equations on $M$ imply that the equations hold on
$N$. Assume that
the forms $E$ and $H$ satisfy Maxwell's equations
(\ref{eq: physical Max BBB})
in $M$ in the classical sense, with traces
(\ref{eq: M restrictions})
that satisfy (\ref{eq: M restrictions same}).
Then $E$ and $H$ are $C^\infty$ smooth near $\gamma$.

Define  1-forms $\tilde E,\tilde H$  and 2-form $\tilde J$
on $N\setminus \Sigma$ by
$\tilde E=(F^{-1})^* E$, $\tilde H=(F^{-1})^*
H$, and $\tilde J=((F^{-1})^* J.$
Then $\tilde E$ satisfies  Maxwell's equations
on
$N\setminus \Sigma$,
\beq\label{eq: physical Max2 B}
\nabla\times
\tilde E = ik \tilde \mu(x)  \tilde H,\quad \nabla\times
           \tilde H
=-ik  \tilde \e(x) \tilde E+\tilde J\quad
\hbox{ in }N\setminus
\Sigma.
\eeq
A simple computation shows that
$\tilde E$, $\tilde H$, $\tilde D=\tilde \e\tilde E$,
and $\tilde B=\tilde \mu\tilde H$
are forms on $N$ with $L^1(N,dx)$ coefficients.
Again, let $\Sigma(t)$ be the
$t$-neighborhood
of
$\Sigma$ in $\tilde g$-metric and
$\gamma(t)$ be
the
$t$-neighborhood
of $\gamma$ in $g$-metric.
Let  $I_t:\p \gamma
(t)\to M$ be the identity
embedding. Denote by
$\nu$ be the
unit normal vector of
$\p\Sigma(t)$ and $\p\gamma(t)$
  in Euclidean metric.

Now, writing $E=E_j(x)dx^j$ on $M$, we see
as
above using $F_t= F\circ
I_t:\p \gamma(t)\to \p \Sigma(t)$, we have
in local coordinates
formula (\ref{tilde I formula}).
Let us next do
computations in the Euclidean
coordinates.
Using
(\ref{rotations}),
the angular direction
$\eta:=\p_\theta$,
and
vertical direction $\zeta=\p_z$,
we see that the matrix
$DF_t^{-1}(x)$
satisfies
\ba
&
&|\eta\,\cdotp
(DF_t^{-1}(x)\eta)|_{\R^3}\leq Ct,\quad x\in
\p\Sigma(t),\\
&
&|\zeta\,\cdotp (DF_t^{-1}(x)\zeta)|_{\R^3}=1
,\quad
x\in
\p\Sigma(t),\\
& &\zeta\,\cdotp (DF_t^{-1}(x)\eta)=0,\quad
x\in
\p\Sigma(t),\\
& &\eta\,\cdotp (DF_t^{-1}(x)\zeta)=0,\quad
x\in
\p\Sigma(t).
\ea
This implies that only angular components of
$\tilde
E$ vanish
on $\Sigma$, and we have

\beq\label{point limits}
& &|\eta\,\cdotp \tilde E|_{\R^3}\leq Ct,\quad x\in
\p\Sigma(t),\\
& &\lim_{t\to 0}\zeta\,\cdotp \tilde
E|_{\p\Sigma(t)\cap N_j} =\tilde
b^e_j,\quad j=1,2 \nonumber
,\\
& &\lim_{t\to
0}\zeta\,\cdotp \tilde
H|_{\p\Sigma(t)\cap N_j}=\tilde
b^h_j,\quad j=1,2 \nonumber
\eeq
where, for $(x^1,x^2,x^3)\in \Sigma\subset
N$, we
denote
\ba
\tilde
b^e_j(x^1,x^2,x^3)=b^e_j(x^3),\quad
\tilde
b^h_j(x^1,x^2,x^3)=b^h_j(x^3),\quad
j=1,2.
\ea

Thus, using
(\ref{eq: physical Max2 B})
we see that for
$\tilde h\in
C^{\infty}_0(\Omega^1N)$
\beq\label{eq: weak N 2 Max
B}
&
&\int_N
((\nabla\times \tilde h)\,\cdotp \tilde E-
           ik
\tilde h
\,\cdotp \tilde \mu(x)\tilde H) \,dV_0(x)\\
& &=\lim_{t\to
0}
\int_{N\setminus \Sigma(t)}
((\nabla\times \tilde h)\,\cdotp \tilde
E-ik \tilde h \,\cdotp \tilde \mu(x)\tilde H)
\,dV_0(x)
\nonumber
\\
& &=-\lim_{t\to 0} \int_{\p \Sigma(t)}
(\nu
\times \tilde E)\,\cdotp \tilde h\,  dS(x), \nonumber
\\
&
&=-\int_{\Sigma}
(\nu \times (\tilde b^e_1-\tilde b^e_2)\zeta)\,\cdotp
\tilde h\,
dS(x) \nonumber\\
\nonumber
& &=0\eeq
where $\nu$ is the Euclidian unit normal
of $\p N_2=\Sigma$.
This shows that Maxwell's equations are satisfied on $N$.
     Observe  that if $\tilde b^e_1\not =\tilde b^e_2$,
there exists a test function $\tilde h$ such that the last integral is nonzero,
precluding the existence of a finite energy solution.
Similar considerations are valid for the  equation
$\nabla \times \tilde H=-ik\tilde \e\tilde E+\tilde J$.

On the other hand, assume that
1-forms $ \tilde E$ and $ \tilde H$
satisfy on $N$ Maxwell's equations (\ref{eq: physical Max tube})
in the finite energy sense.
Then, as $E$ and $H$ are forms with $L^2(M)$-valued
coefficients that satisfy Maxwell's equations in
$M_1\setminus \gamma_1$ and $M_2\setminus \gamma_2$,
we see that they have to satisfy
Maxwell's equations in
$M_1$ and $M_2$, and thus they are $C^\infty$-smooth forms
near $\gamma_1$ and $\gamma_2$.
As $\tilde E$ and $\tilde H$ are finite energy solutions on $N$,
the above arguments show that
$ b^e_1= b^e_2$
and   $b^h_1= b^h_2$. This finishes the proof of
Theorem \ref{single coating with Maxwell}.
\hfill\proofbox

{\bf Remark.} If $E$, $H$, and $J$ on $M$ are solutions of Maxwell's equations
as in Proposition \ref{single coating with Maxwell} (1)
such that conditions (114) are not satisfied,
then the proof of Proposition 7.1
shows that the fields $\tilde E=F_*E$, $\tilde H=F_*H$, and $\tilde J=F_*J$
extend to forms with $L^1(N,dx)$ coefficients that satisfy
\beq\label{eq: new equations}
& &\nabla\times \tilde E =ik \tilde B+\tilde K_{new}\quad \hbox{on }N ,\\
& &\nabla\times \tilde H =-ik \tilde D+\tilde J+\tilde J_{new}\quad\hbox{on }N
\nonumber
\eeq
in the sense of distributions. Here, $\tilde B=\tilde \mu\tilde H$
and $\tilde D=\tilde \e\tilde E$ are 2-forms with measurable coefficients
and
$\tilde K_{new}=s_e\delta_{\Sigma}$
and $\tilde J_{new}=s_h\delta_{\Sigma}$
where $\delta_\Sigma$ is a measure supported on $\Sigma$ and $s_e$ and
$s_h$ are smooth 2-forms.

Similarly, if
$E$, $H$, and $J$ on $M$ are solutions of Maxwell's equations
as in Theorem \ref{single coating with Maxwell} (2) with non-vanishing
Cauchy data (\ref{eq: physical Max single M3 new}), we see that
that $\tilde E$, $\tilde H$, and $\tilde J$ on $N$ satisfy equations (\ref
{eq: new equations}) with distributional sources
$\tilde K_{new}$ and $\tilde J_{new}$ defined as above.

\section{Cloaking a cylinder with
the SHS  boundary condition}\label{cylsec-shs}

Next, we consider  $N_2$ as an obstacle,
while the domain  $N_1$ is equipped with a metric corresponding to
the single coating.
Motivated by the conditions at $\Sigma$ in the previous section,
we impose the soft-and-hard  boundary condition on
the boundary of the obstacle. To this end, let us  give
still one more definition of weak solutions, appropriate for this construction.
We consider only solutions on the set $N_1$;
nevertheless, we continue to denote $\p N=\p N_1\setminus \Sigma$.

\begin{definition}\label{SHS-def}
Let $(M_1,N_1,F,\gamma_1,\Sigma,g_1)$ be a single coating construction.
We say that the 1-forms $ \tilde E$ and $ \tilde H$ are \emph{finite 
energy solutions } of  Maxwell's equations on $N_1$
with the soft-and-hard (SHS) boundary conditions
on $\Sigma$,
\beq\label{eq: physical Max B Ob}
& &\nabla\times \tilde E = ik\tilde \mu(x)  \tilde H,
\quad \nabla\times \tilde H =-ik  \tilde \e(x) \tilde E+\tilde J\quad
\hbox{ on }N_1,\\
\label{special electric and magnetic condition 2 A Ob}
& &\eta\,\cdotp \tilde E|_{\Sigma }=0,\quad
\eta\,\cdotp \tilde H|_{\Sigma}=0,
\\
\nonumber
& &\nu \times  \tilde E|_{\p N}=f,
\eeq
if $ \tilde E$ and $ \tilde H$ are 1-forms on $N_1$ and $\tilde \e 
\tilde E$ and $\tilde \mu \tilde H$ are 2-forms with
measurable coefficients
satisfying
\beq\label{eq: Max norm 3Ob}
\|\tilde E\|_{L^2(N_1,|\tilde g|^{1/2}dV_0)}^2=
\int_{N_1} \tilde \e^{ij}\, \tilde E_j\, \overline{\tilde E_k}
\,dV_0(x)<\infty,\\
\|\tilde H\|_{L^2(N_1,|\tilde g|^{1/2}dV_0)}^2=\int_{N_1}
\tilde \mu^{ij}\, \tilde H_j\, \overline{\tilde H_k}
\,dV_0(x)<\infty;
\eeq
Maxwell's equation are valid in
the classical sense in a
neighborhood $U$ of $\p N$:
\ba
& &\nabla\times  \tilde E = ik  \tilde \mu(x)   \tilde H,\quad
\nabla\times \tilde  H =-ik  \e(x)\tilde  E+\tilde J\quad
\hbox{in }U,\\
& &\nu\times \tilde E|_{\p N}=f;
\nonumber
\ea
and
finally,
\ba
& &\int_{N_1} ((\nabla\times \tilde h)\,\cdotp \tilde E-
            ik \tilde h \,\cdotp \tilde \mu(x)\tilde H)
\,dV_0(x)=0,\\
& &\int_N ((\nabla\times \tilde e)\,\cdotp \tilde H+
\tilde e \,\cdotp ( ik\tilde \e(x)\tilde E-\tilde J))
\,dV_0(x)=0,
\ea
for all $\tilde e,\tilde h\in C^\infty_0(\Omega^1N_1)$
satisfying
\beq\label{special electric and magnetic condition 2 A1}
\eta\,\cdotp \tilde e|_{\Sigma}=0,\quad \eta\,\cdotp \tilde h|_{\Sigma}=0,
\eeq
where $\eta=\p_\theta$ is the angular vector field that is tangential to
$\Sigma$.
\end{definition}

We have the following invisibility result.

{In this section  $(M_1,N_1,F,\gamma_1,\Sigma)$ is
a coating configuration corresponding to single coating of
a cylindrical obstacle $B_2(0,1)\times \R$.}

\begin{theorem}\label{single coating with Maxwell obstacle}
Let $E$ and $H$ be 1-forms
with measurable coefficients
on $M_1\setminus \gamma_1$ and $\tilde E$ and $\tilde H$ be
           1-forms
with measurable coefficients on $N_1$
such that $E=F^*\tilde E$, $H=F^*\tilde H$.
Let $J$ and $\tilde J$ be  2-forms
with smooth coefficients on $M_1\setminus \gamma_1$ and
           $N_1\setminus \Sigma$, that are supported away
from $\gamma_1$ and $\Sigma$.

Then the following are equivalent:
\begin{enumerate}

\item
On $N_1$, the
1-forms $ \tilde E$ and $ \tilde H$ satisfy Maxwell's equations
(\ref{eq: physical Max B Ob}) with SHS boundary conditions
(\ref{special electric and magnetic condition 2 A Ob})
in the sense of Definition \ref{SHS-def}.

\item
On $M_1$, the forms $E$ and $H$ are classical solutions of Maxwell's equations,
\beq\label{eq: physical Max BBB2}
& &\nabla\times     E = ik    \mu(x)   H,\quad
\hbox{in }M_1\\
\nonumber
& &\nabla\times H =-ik  \e(x)   E + J,\quad
\hbox{in }M_1,
\\
& &\nu\times
E|_{\p M_1}=f. \nonumber
\eeq

\end{enumerate}
\end{theorem}

\noindent
{\bf Proof.} First, assume that the forms $E$ and $H$ satisfy
Maxwell's equations
(\ref{eq: physical Max BBB2}) in $M_1$. Then $E$ satisfies
identities (\ref{point limits}).
Considerations similar to those yielding formula (\ref{eq: weak N 2
Max B}) imply that
\beq\label{eq: weak N 2 Max B Ob}
&
&\int_{N_1}
((\nabla\times \tilde h)\,\cdotp \tilde E-
            ik \tilde h
\,\cdotp \tilde \mu(x)\tilde H) \,dV_0(x)\\
& &=\lim_{t\to 0}
\int_{{N_1}\setminus \Sigma(t)}
((\nabla\times \tilde h)\,\cdotp \tilde
E-
            ik \tilde h \,\cdotp \tilde \mu(x)\tilde H)
\,dV_0(x)
\nonumber
\\
& &=-\lim_{t\to 0} \int_{\p \Sigma(t)\cap {N_1}}
(\nu
\times \tilde E)\,\cdotp \tilde h\,  dS(x), \nonumber
\\
& &=-\lim_{t\to 0} \int_{\p \Sigma(t)\cap {N_1}}
(\nu
\times
((\eta\,\cdotp \tilde E)\eta +(\zeta\,\cdotp \tilde E)\zeta
)\,\cdotp
\tilde h\,  dS(x), \nonumber
\\
& \nonumber
&=0
\eeq
for atest
function $\tilde h$ satisfying
(\ref{special electric and magnetic
condition 2 A1}).

Similar analysis for $\tilde H$ shows
that
1-forms $ \tilde E$ and $ \tilde H$ satisfy Maxwell's
equations
with SHS boundary conditions
in the sense of Definition
\ref{SHS-def}.

Next, we show that
equations on
$N_1$ imply
equations on $M_1$.
Assume that 1-forms $ \tilde E$ and $ \tilde H$
satisfy Maxwell's equations
with SHS boundary conditions, and
internal current $\tilde J$,
in the sense of Definition
\ref{SHS-def}. Then $E$ and $H$ are
classical solutions of Maxwell's
equation in $M_1\setminus \gamma_1$.
Let
$U\subset M_1$ be a
neighborhood
of $\gamma_1$ and
$W=F(U\setminus \gamma_1)\cup \Sigma$
be
a neighborhood of $\Sigma$ in $N_1$
such that
$\hbox{supp}\,(\tilde J)\cap W=
\emptyset$.
Then we have
\ba
\tilde
\e^{jk}\tilde E_j\overline {\tilde E_k} \in
L^1(W,\,dV_0(x)),\quad
  \tilde \mu^{jk}\tilde H_j\overline {\tilde H_k}
\in
L^1(W,\,dV_0(x)).
\ea
Define  $E=F^*\tilde E$,
     $H=F^*\tilde H$
and
$J=F^*\tilde J$ on $M_1\setminus \gamma_1$.
Again, we
see
that
$E,$ $H$, and $J$ satisfy Maxwell's equations on $U\setminus
\gamma$,
and as above we see that $E$ and $H$ have
measurable
extensions on $\gamma$,
$E^e,H^e\in L^2(U,\,dV_0(x))$
, such
that
$\nabla\times E^e -ik \mu(x) H^e$ and
$\nabla\times H^e +ik
\e(x) E^e$ are distributions in
$H^{-1}(U,dV_0)$ supported on
$\gamma_1$.
As before, we see obtain
\ba
        \nabla\times E^e -ik
\mu(x) H^e=0,\quad
          \nabla\times H^e +ik \e(x) E^e=0\quad
\hbox{in } U.
\ea
This shows that $E$ and $H$ are classical solutions
of Maxwell's equations
on $M_1$.\hfill\proofbox
\medskip

Similar
analysis can be done in the case
when we have a physical surface
$\Sigma=S^1\times \R$ dividing  $\R^3$
into two regions, having the
SHS boundary conditions
on  both sides, and we define the material
parameters according
to  double coating construction, i.e., on both
sides of
the surface.

\section{Appendix: Single and double coating
for arbitrary domains
and metrics}\label{appendix}

The
constructions of \S2 and the results that follow easily extend
to
general domains and metrics.
Let us assume that
$\Omega\subset\R^3$ now is an arbitrary domain with smooth
boundary,
equipped with an
arbitrary smooth Riemannian metric, $g=g_{ij}(x)$.
This defines the
Laplace operator
$\Delta_g$ with, say Dirichlet
boundary condition, cf. Remark 3.6.
Choose a point ${O} \in \Omega$
to be blown up,  and
         assume that the  injectivity radius

of $(\Omega, g)$ at ${O}$ is larger than $3a$ for some $ a>0$.
Let
$B(O,r)$ denote a metric ball of $(M,g)$
    with center $O$ and
radius $r$.
        Introduce Riemannian normal coordinates in
$B({O},3a)
\subset \Omega:$
         $$
         x=(x^1,x^2,x^3)
\to(\tau, \omega), \tau >0, \omega \in \stwo\subset
T_{O}\Omega,
  $$
so that $x =\exp_{{O}}(\tau \omega)$.
Let $f(\tau):[0,3a] \to
[a,3a]$ be a smooth strictly increasing function
coinciding
with $\tau/2+a$ near $\tau =0$ and with $\tau$ for
$\tau >2a$.

Define, in these coordinates,
         $$
         F: B({O},3a)
\setminus \{{O}\} \to
         B({O},3a) \setminus
         B({O},a),
  \quad
         (\tau, \omega) \to (f(\tau), \omega).
         $$
  We extend $F$
by the identity to $\Omega \setminus B({O},3a)$
and obtain a
diffeomorphism
         $$
         F_1: \Omega \setminus
\{{O}\} \to N_1=\Omega \setminus
B({O},a).
         $$
Consider the metric $\tilde g = {F_1}_* g$ in $N_1$. Observe that
surfaces lying at distance $\tau$ from $\p  B({O},a)$ with respect
to the metric
         $\tilde g$ coincide with surfaces lying at
distance $f(\tau)-a$ from
         $\p B({O},a)$ with respect to the
metric $g$.
Therefore,  the directions normal to these
surfaces are  the
same with respect to the metrics $g$ and
$\tilde g$.
In particular, the direction of these normals, in
the metric $\tilde g$,
is transversal to  $\p  B({O},a)$.
Thus,
equations (\ref{propertyextra}) remain valid if we use
$\tau - a$
instead of $r-1$.
Similarly, we again have the
estimate $|\tilde g|^{1/2} \leq C_1(\tau-a)^2$.

One may also extend
the double coating construction as follows.
Let  $(D, g_D)$
be a  compact Riemannian manifold without boundary, and
choose a
point
$NP \in D$.
Using Riemannian normal coordinates centered
at $NP$,  introduce,
similar to the
above, a diffeomorphism
  $$
         F_2: D \setminus \{NP\} \to N_2=D \setminus
\overline{B}(NP,b),
         $$
         where we assume that $3b $ is
smaller than injectivity radius of $D$.
         Pulling back the
metric $g_D$, we get a
         metric $\tilde{g_D}$ on $D \setminus
\overline{B}(NP,b)$ with the
same properties near $\p B(NP,b)$ as
$\tilde g$ has near $\p
B({O},a)$.

Observe that, as we
are inside
the injectivity radii,
$\p B({O},a)$ and $\p B(NP,b)$ are both
diffeomorphic to
$\stwo$,  with
diffeomorphisms given by
  $exp_{{O}}(a\omega)$ and $\exp_{NP}(b\omega)$.  Thus,
          $\p
B({O},a)$ and $\p B(NP,b)$  are diffeomorphic to each
         other.
  Gluing these boundaries, we obtain a smooth manifold $N=N_1
\cup
N_2\cup\Sigma$
         with a  Riemannian metric singular on
$\Sigma$
          which, as one approaches $\Sigma$, satisfies
conditions
(\ref{propertyextra}).
          This makes it possible to
carry out all   of the preceding analysis
         for the double
coating.

         Note that if $D$ is diffeomorphic to $S^3$ (as
earlier), then
$N$ is diffeomorphic to $\Omega\simeq M_1$. If
however $D$ has a
non-trivial topology, $N$ may have topology
different from that of
$\Omega$. However, due to the full-wave
invisibility, one is unable
to observe this change of  topology
from observations made at $\p \Omega$. Note that this is
in contrast to the uniqueness result that holds for $C^\omega$ Riemannian
manifolds \cite{LTU}.

Similar
generalizations of the single coating construction are
possible when $\p D$ is diffeomorphic to $\stwo$.

\bibliographystyle{amsalpha}

\bibliographystyle{amsalpha}

\end{document}